\documentclass[10pt]{article}

\usepackage{a4wide}
\usepackage{amssymb}
\usepackage{amsfonts}
\usepackage{amsmath}
\input xy
\xyoption{arrow} \xyoption{matrix}

\date{}

\newtheorem{proposition}{Proposition}[section]
\newtheorem{theorem}[proposition]{Theorem}
\newtheorem{lemma}[proposition]{Lemma}

\newtheorem{corollary}[proposition]{Corollary}

\def\Hom{{\rm Hom}}
\def\der{\partial }

\def\nFM0{{\nu }_{F,M_0}}
\def\nFN0{{\nu }_{F,N_0}}
\def\nGN0{{\nu }_{G,N_0}}

\def\N0{ {\bf N}_0 }

\def\g{\gamma}

\def\ra{\rightarrow}

\def\Xpm{X^{\pm }}

\def\s{\sigma}
\def\Z{\mathbb{Z}}

\def\l1{{\lambda}_1}

\def\a{\alpha}
\def\a0{ {\alpha }_0}
\def\a1{ {\alpha }_1}

\def\l{\lambda}

\def\nFGM0{{\nu }_{F,G,M_0}}

\def\nFN0{{\nu}_{F,N_0}}

\def\sm{{\sigma}^m}

\def\sm1{{\sigma}^{-1}}

\def\smtp1{{\sigma}^{-t+1}}

\def\S1{S^{-1}}

\def\Xpm1{X^{\pm 1}_1}

\def\sPM1{{\sigma }^{\pm 1}}
\def\sMP1{{\sigma }^{\mp 1 }}

\def\d{\delta}

\def\di{{\rm d.ind}}

\def\L{\Lambda}

\def\CD{{\cal D}}

\def\Ytm1{Y^{t-1}}
\def\Yim1{Y^{i-1}}

\def\CM{{\cal M}}
\def\CN{{\cal N}}

\def\CF{{\cal F}}
\def\CG{{\cal G}}
\def\CH{{\cal H}}

\def\Aut{{\rm Aut}}

\def\Der{{\rm Der }}
\def\ad{{\rm ad }}
\def\dim{{\rm dim }}

\def\ker{ {\rm ker } }

\def\CJ{ {\cal J}}

\def\gcd{ {\rm gcd } }
\def\D{ \Delta }

\def\SL2Z{ {\rm SL}_2({\bf Z}) }

\def\th{ \theta }

\def\Gp1{ G^{1 , 1 } }
\def\P11{ P^{-1 , 1 } }
\def\Pp1{ P^{1 , 1 } }

\def\Supp{{\rm Supp}}

\def\th{\theta}

\def\nCLsr{{}^\nu\kern-2pt {\cal L}^{\sigma , \rho  }}
\def\nP{{}^\nu \kern-2pt P}
\def\nL{{}^\nu\kern-2pt L}
\def\nLL{{}^\nu\kern-2pt \Lambda}
\def\nPsr{{}^\nu\kern-2pt P^{\sigma , \rho  }}
\def\nLsr{{}^\nu\kern-2pt L^{\sigma , \rho  }}
\def\nuCL{{}^\nu\kern-2pt  {\cal L}}
\def\nCLsr{{}^\nu\kern-2pt {\cal L}^{\sigma , \rho  }}
\def\nCL1m{{}^\nu\kern-2pt {\cal L}^{-1 , 1  }}
\def\x1nu{x^\frac{1}{\nu}}
\def\xm1nu{x^{-\frac{1}{\nu}}}

\def\CN{{\cal N}}
\def\ra{\rightarrow }

\def\CB{{\cal B}}

\def\CI{{\cal I}}

\def\CH{ {\cal H}}

\def\nAM0{{\nu }_{{\cal A},M_0}}
\def\nAN0{{\nu }_{{\cal A},N_0}}

\def\End{ {\rm End }}
\def\Der{ {\rm Der }}
\def\CJ{ {\cal J }}

\def\det{ {\rm det }}
\def\ad{ {\rm ad }}

\def\ga{\mathfrak{a}}

\def\GL{{\rm GL}}
\def\SL{{\rm SL}}

\def\Hom{{\rm Hom}}

\def\di!{\frac{\der^i}{i!}}
\def\dik!{\frac{\der^k_i}{k!}}

\def\id{{\rm id}}

\def\N{\mathbb{N}}

\def\0{\overline{0}}
\def\1{\overline{1}}

\def\Ln1{\L_{n,\overline{1}}}

\def\a1{a_{\overline{1}}}

\def\S{\Sigma}

\def\grad{{\rm grad}}

\def\vn1{\overrightarrow{n-1}}

\def\Sh{{\rm Sh}}

\def\Inn{{\rm Inn}}

\def\mJ{\mathbb{J}}
\def\mI{\mathbb{I}}

\def\mT{\mathbb{T}}

\def\mG{\mathbb{G}}

\def\K1{{\rm K}_1}

\def\hmI1{\widehat{\mI_1}}
\def\tmI1{\widetilde{\mI_1}}
\def\tmJ1{\widetilde{\mJ_1}}
\def\hB1{\widehat{B_1}}
\def\hCB1{\widehat{\CB_1}}

\def\Fix{{\rm Fix}}

\def\mJ{\mathbb{J}}

\def\AutKalg{ {\rm Aut_{K-{\rm alg}}}}

\def\mmW1{\mathbb{W}_1}

\def\divn0{\mathfrak{div}_n^0}
\def\div0mu{\mathfrak{div}_{n, [\mu ]}^0}
\def\din0{\mathfrak{di}_n^0}
\def\ivn0{\mathfrak{iv}_n^0}
\def\divnc{\mathfrak{div}_n^c}
\def\divv{{\rm div}}
\def\mfGn{\mathbf{G}_n}
\def\mfGnc{\mathbf{G}_n^c}
\def\mh{\mathfrak{h}}
\begin{document}

\author{V. V. \  Bavula   
}

\title{The groups of automorphisms of the Lie algebras of polynomial vector fields with zero or constant divergence}

\maketitle

\begin{abstract}
Let $P_n=K[x_1, \ldots , x_n]$ be a polynomial algebra over a field $K$ of characteristic zero and $\divn0$ (respectively, $\divnc$) be the Lie algebra of derivations of $P_n$ with zero (respectively, constant) divergence.
We prove that $\Aut_{{\rm Lie}}(\divn0 )\simeq \Aut_{K-{\rm alg}}(P_n)$ ($n\geq 2$) and $\Aut_{{\rm Lie}}(\divnc )\simeq \Aut_{K-{\rm alg}}(P_n)$. The Lie algebra $\divnc$ is  a maximal Lie subalgebra of $\Der_K(P_n)$. Minimal  finite sets of generators are found  for the Lie algebras $\divn0$ and $\divnc$.

$\noindent $



\end{abstract}


\section{Introduction}

In this paper, module means a left module, $K$ is a
field of characteristic zero and  $K^*$ is its group of units, and the following notation is fixed:
\begin{itemize}
\item $P_n:= K[x_1, \ldots , x_n]=\bigoplus_{\alpha \in \N^n}
Kx^{\alpha}$ is a polynomial algebra over $K$ where
$x^{\alpha}:=x_1^{\alpha_1}\cdots x_n^{\alpha_n}$ and $Q_n:= K(x_1,\ldots , x_n)$ is its field of fractions,
 \item $G_n:=\AutKalg (P_n)$ is the group of automorphisms of the polynomial algebra $P_n$,
 \item $\der_1:=\frac{\der}{\der x_1}, \ldots , \der_n:=\frac{\der}{\der
x_n}$ are the partial derivatives ($K$-linear derivations) of
$P_n$,
\item    $D_n:=\Der_K(P_n) =\bigoplus_{i=1}^nP_n\der_i$ is the Lie
algebra of $K$-derivations of $P_n$ where $[\der , \d ]:= \der \d -\d \der $,
\item $\mG_n:=\Aut_{{\rm Lie}}(D_n)$ is the group of automorphisms of the Lie algebra $D_n$,
\item  $\d_1:=\ad (\der_1), \ldots , \d_n:=\ad (\der_n)$ are the inner derivations of the Lie algebra $D_n$ determined by  $\der_1, \ldots , \der_n$ (where $\ad (a)(b):=[a,b]$),
 \item $\CD_n:=\bigoplus_{i=1}^n K\der_i$,
 \item $\CH_n :=\bigoplus_{i=1}^n KH_i$ where $H_1:=x_1\der_1, \ldots , H_n:=x_n\der_n$,
     \item $D_n':=\bigoplus_{i=1}^nP_nH_i= \bigoplus_{\alpha\in \N^n} x^\alpha \CH_n$,
  \item $\mh :=\bigoplus_{i=1}^n Kh_i$ where $h_1:=\der_1x_1, \ldots, h_n:=\der_nx_n\in \End_K(P_n)$,

 \item for a derivation $\der = \sum_{i=1}^n a_i\der_i\in D_n$, $\divv (\der ) := \sum_{i=1}^n\frac{\der  a_i}{\der x_i}$ is the {\em divergence}  of $\der$,
\item $\divn0 :=\{ \der \in D_n \, | \, \divv (\der ) =0\}$ is the Lie algebra  of polynomial vector fields (derivations) with zero divergence,
 \item $\mfGn :=\Aut_{{\rm Lie}}(\divn0 )$,
  \item $\CH_n' :=\bigoplus_{i=1}^{n-1} KH_{i,i+1}$ where $H_{ij} :=H_i-H_j$ for $i\neq j$,
 \item $\divnc :=\{ \der \in D_n \, | \, \divv (\der ) \in K\}$ is the Lie algebra  of polynomial vector fields (derivations) with constant  divergence,

 \item $\mfGnc :=\Aut_{{\rm Lie}}(\divnc )$,
  \item  $A_n:= K \langle x_1, \ldots
, x_n , \der_1, \ldots , \der_n\rangle  =\bigoplus_{\alpha , \beta
\in \N^n} Kx^\alpha \der^\beta$  is  the $n$'th {\em Weyl
algebra},
    \end{itemize}

{\bf The groups of automorphisms of the Lie algebras $\divn0$ and $\divnc$}.
The aim of the paper is to prove the following two theorems.

\begin{theorem}\label{16Mar13}
$\mfGn = \begin{cases}
G_1/\Sh_1\simeq K^*& \text{if }n=1,\\
G_n& \text{if }n\geq 2.\\
\end{cases}$
\end{theorem}
{\em Structure of the proof}. The case $n=1$ is trivial (see Section \ref{P16AAA} where the group $\Sh_1$ is defined in (\ref{Shndef})). So, let $n\geq 2$.

 $\noindent $

 (i) $G_n\subseteq \mfGn$ via the group monomorphism (Lemma \ref{db11Mar13}.(3))
$$G_n\ra \mfGn, \;\;  \s \mapsto \s : \der \mapsto \s (\der ):=\s \der \s^{-1}.$$


(ii) Let $\s \in \mfGn$. Then $\der_1':=\s (\der_1), \ldots , \der_n':=\s (\der_n)$ are commuting, locally nilpotent derivations of the polynomial algebra $P_n$ (Lemma \ref{dc13Mar13}.(1)).

$\noindent $

(iii) $\bigcap_{i=1}^n \ker_{P_n}(\der_i') = K$  (Lemma \ref{dc13Mar13}.(2)).

$\noindent $

(iv)  There exists a polynomial automorphism $\tau \in G_n$ such that $\tau \s \in \Fix_{\mfGn}(\der_1, \ldots , \der_n)$ (Corollary \ref{db13Mar13}).

$\noindent $

(v) $\Fix_{\mfGn}(\der_1, \ldots , \der_n)=\Sh_n$ (Proposition \ref{dB11Mar13}.(3)) where
\begin{equation}\label{Shndef}
\Sh_n:=\{ s_\l \in G_n\, | \, s_\l (x_1)=x_1+\l_1, \ldots , s_\l (x_n) = x_n+\l_n\}
\end{equation}
is the {\em shift group} of automorphisms of the polynomial algebra $P_n$ and $\l = (\l_1, \ldots , \l_n)\in K^n$.

$\noindent $

(vi) By (iv) and (v), $\s \in G_n$, i.e. $\mfGn = G_n$.   $\Box $


\begin{theorem}\label{A16Mar13}
$\mfGnc = G_n$.
\end{theorem}
{\em Structure of the proof}. The case $n=1$ is trivial (see Section \ref{P16AAA}). So, let $n\geq 2$.

 $\noindent $

  (i) $G_n\subseteq \mfGnc$ via the group monomorphism (Lemma \ref{db11Mar13}.(4))
$$G_n\ra \mfGnc, \;\;  \s \mapsto \s : \der \mapsto \s (\der ):=\s \der \s^{-1}.$$


(ii) $\divn0 = [\divnc , \divnc ]\;\;$ (Lemma \ref{a16Mar13}).

$\noindent $

(iii) The short exact sequence of group homomorphisms
$$ 1\ra F:=\Fix_{\mfGnc}(\divn0 ) \ra \mathbf{G}_n^c
\stackrel{{\rm res}}{\ra}\mfGn\ra 1$$
is exact  (by (i) and Theorem \ref{16Mar13}) where res $:\s \mapsto \s|_{\divn0}$ is the restriction map, see (ii).

 $\noindent $

 (iv)  Since $\mfGn = G_n$ (Theorem \ref{16Mar13}) and $G_n\subseteq \mfGnc$ (by (i)), the short exact sequence splits
\begin{equation}\label{G=GF}
 \mfGnc \simeq G_n\ltimes F.
\end{equation}


(v) $ F=\{ e\}$ (Lemma \ref{a18Mar13}). Therefore, $\mfGnc = G_n$.   $\Box $

$\noindent $




Theorem \ref{16Mar13} was announced  in \cite{Rudakov-1986}  where a  sketch of the proof is given based on a study of certain Lie subalgebras of $\divn0$ of finite codimension. Our proof is based on completely different ideas. The groups of automorphisms of infinite dimensional Lie algebras were considered in [2]-[10].

$\noindent $

A subalgebra $\CM$ of a Lie algebra $\CG$ is called a {\em maximal} Lie subalgebra if $\CM \neq \CG$ and  $\CG$ is the only   Lie subalgebra of $\CG$ properly containing $\CM$.
\begin{itemize}
\item {\rm (Proposition \ref{a2Apr13})} {\em For $n\geq 2$,  $ \divnc$ is  a maximal Lie subalgebra of} $D_n$ which is also a $\mG_n$-invariant/$G_n$-invariant Lie subalgebra.

\item {\rm (Proposition \ref{c2Apr13})} {\em For $n\geq 2$, the $G_n$-module $D_n/ \divnc$ is simple and infinite dimensional with} $\End_{G_n}(D_n/\divnc )\simeq K$.
    \end{itemize}

\begin{theorem}\label{6Oct13}
For $n\geq 2$, the set of  elements $ \der_1, x_2^2\der_1,  x_3^2\der_1,\ldots ,  x_n^2\der_1,
 x_1^2\der_2,  x_1^2\der_3, \ldots,  x_1^2\der_n$ is a minimal set of generators for the Lie algebra $\divn0$.
\end{theorem}

\begin{theorem}\label{A6Oct13}
For $n\geq 2$, the set of elements in Theorem \ref{6Oct13} together with  $H_1$
 is a minimal set of generators for the Lie algebra $\divnc$.
\end{theorem}

\section{Proof of Theorems \ref{16Mar13}  and \ref{A16Mar13}}\label{P16AAA}

In this section, proofs of Theorems \ref{16Mar13}  and \ref{A16Mar13} are given. In the first part of the section some useful results are proved that are used throughout the paper. The second part of the  section can be seen as  proofs of Theorem \ref{16Mar13}  and \ref{A16Mar13}. The proofs are  split into several statements that reflect `Structure of the proofs of Theorems \ref{16Mar13} and  and \ref{A16Mar13}'  given in the Introduction. As we have seen in the Introduction, Theorem \ref{16Mar13} is the key point in the proof of Theorem \ref{A16Mar13}.

$\noindent $

{\bf The Lie algebra $D_n$ is $\Z^n$-graded}. The Lie algebra
\begin{equation}\label{xadbd}
D_n =\bigoplus_{\alpha\in \N^n} \bigoplus_{i=1}^n Kx^\alpha \der_i
\end{equation}

 is a $\Z^n$-graded Lie algebra
$$D_n = \bigoplus_{\beta\in \Z^n} D_{n , \beta}\;\; {\rm where}\;\; D_{n,\beta}=\bigoplus_{\alpha -e_i=\beta}Kx^\alpha \der_i,$$
i.e. $[D_{n,\alpha} , D_{n,\beta }] \subseteq D_{n, \alpha +\beta}$ for all $\alpha ,\beta \in \N^n$ where $e_1:=(1, 0 , \ldots , 0), \ldots , e_n:=(0, \ldots , 0 , 1)$ is the canonical free basis for the free abelian group $\Z^n$. This follows from the commutation relations

\begin{equation}\label{xadbd1}
[x^\alpha\der_i, x^\beta \der_j]= \beta_i x^{\alpha+\beta - e_i} \der_j-\alpha_j x^{\alpha + \beta - e_j} \der_i.
\end{equation}
Clearly, for all $i,j=1, \ldots , n$ and $\alpha \in \N^n$,
\begin{equation}\label{xadbd2}
[H_j, x^\alpha \der_i]=\begin{cases}
\alpha_j x^{\alpha} \der_i & \text{if }j\neq i ,\\
(\alpha_i-1)x^{\alpha} \der_i& \text{if }j=i, \\
\end{cases}
\end{equation}
\begin{equation}\label{xadbd3}
[\der_j, x^\alpha \der_i]=\alpha_j x^{\alpha -e_j} \der_i.
\end{equation}
The {\em support} $\Supp (D_n):=\{ \beta \in \Z^n\, | \, D_{n,\beta}\neq 0\}$ is a submonoid of $\Z^n$. Let us find the support $\Supp (D_n)$, the graded components $D_{n,\beta}$ and their dimensions $\dim_K\, D_{n,\beta}$. For each $i=1, \ldots , n$, let $\N^{n,i}:=\{ \alpha \in \N^n \, | \, \alpha_i=0\}$ and $P_n^{\der_i}:=\ker_{P_n}(\der_i)$. It follows from the decompositions
$P_n = P_n^{\der_i}\oplus P_nx_i$ for $i=1, \ldots , n$ that
\begin{equation}\label{Dnb2}
D_n = \bigoplus_{i=1}^n (P_n^{\der_i}\oplus P_nx_i)\der_i =\bigoplus_{i=1}^n P_n^{\der_i}\der_i \oplus \bigoplus_{i=1}^nP_nH_i=\bigoplus_{i=1}^n P_n^{\der_i}\der_i \oplus \bigoplus_{\alpha \in \N^n} x^\alpha \CH_n,
\end{equation}
Therefore, any derivation $\der = \sum_{i=1}^n a_i \der_i\in D_n$ is the unique sum (where $a_i= b_ix_i +c_i $,  $ b_i\in P_n$ and $c_i\in P_n^{\der_i}$)
\begin{equation}\label{dus}
\der = \sum_{i=1}^n b_iH_i + \sum_{i=1}^n c_i\der_i.
\end{equation}
Hence,
\begin{equation}\label{Dnb}
\Supp (D_n) =\coprod_{i=1}^n (\N^{n,i}-e_i) \coprod\N^n.
\end{equation}
\begin{equation}\label{Dnb1}
D_{n,\beta} =\begin{cases}
Kx^\alpha\der_i& \text{if }\beta = \alpha - e_i\in \N^{n,i}-e_i,\\
x^\beta \CH_n& \text{if }\beta \in \N^n.
\end{cases}
\end{equation}
$$\dim_K\, D_{n,\beta} =\begin{cases}
1& \text{if }\beta = \alpha - e_i\in \N^{n,i}-e_i,\\
n& \text{if }\beta \in \N^n.
\end{cases}$$

Let $\CG$ be an abelian  Lie algebra and $\CG^*:= \Hom_K(\CG , K)$. A $\CG$-module $M$ is called a {\em weight module}
if
 $$M= \bigoplus_{\l \in \CG^*}M_\l \;\; {\rm where}\;\; M_\l :=\{m\in M\, | \, gm=\l (g) m \;\; {\rm for \; all}\;\; g\in \CG\}.$$
The set $W (M):=\{ \l \in \CG^*\, | \, M_\l\neq 0\}$ is  called the set of  {\em weights}  of $M$.

$\noindent $

{\bf The direct sum $\divn0 = \din0 \oplus \ivn0$}. Recall that $D_n'=\bigoplus_{\alpha \in \N^n} x^\alpha \CH_n$. By (\ref{Dnb2}),
\begin{equation}\label{dus2}
\divn0 = \din0\oplus \ivn0   \;\; {\rm where}\;\; \din0 := \divn0 \cap D_n'\;\; {\rm  and}\;\;  \ivn0 := \bigoplus_{i=1}^n P_n^{\der_i} \der_i.
\end{equation}
We will see that $\din0$ is a Lie subalgebra of $\divn0$ but $\ivn0$ is not for $n\geq 2$.
 Clearly, $\mathfrak{di}_1^0=0$ and $ \mathfrak{div}_1^0 = \mathfrak{iv}_1^0=K\der_1$. There are  inclusions
 \begin{eqnarray*}
 \mathfrak{div}_1^0&\subset &  \mathfrak{div}_2^0 \subset \cdots \subset \divn0\subset \cdots ,\\
 \mathfrak{di}_1^0&\subset &  \mathfrak{di}_2^0 \subset \cdots \subset \din0\subset \cdots ,\\
 \mathfrak{iv}_1^0&\subset &  \mathfrak{iv}_2^0 \subset \cdots \subset \ivn0\subset \cdots .\\
\end{eqnarray*}

The $K$-linear maps $h_i=\der_1x_1, \ldots , h_n=\der_nx_n\in \End_K(P_n)$ are bijections since for all $\alpha \in \N^n$ and $i=1, \ldots , n$,
\begin{equation}\label{dus3}
h_i(x^\alpha ) = (\alpha_i+1) x^\alpha .
\end{equation}
The elements $h_1, \ldots , h_n$ commute,  the polynomial algebra $P_n$ is a weight $\mh$-module where  $\mh :=\oplus_{i=1}^n Kh_i$ is an abelian Lie subalgebra of the Lie algebra $\End_K(P_n)$ (where $[f,g]:=fg-gf$) and the set $W(P_n)$ of weights of the $\mh$-module $P_n$ is equal to $(1, \ldots , 1)+\N^n$, i.e. $W(P_n) = \{ \l =(\l_1, \ldots , \l_n) \, | \, \l \in (1, \ldots , 1)+\N^n\}$ where $\l (h_i) = \l_i$ for all $i$. For each derivation $\der = \sum_{i=1}^n a_iH_i\in D_n'$,
\begin{equation}\label{dus4}
\divv (\der ) = \sum_{i=1}^n h_i(a_i).
\end{equation}
{\bf $K$-bases for $\divn0$ and $\divnc$}.
For each pair $i\neq j$, the $K$-linear map
\begin{equation}\label{dus5}
\phi_{ij}:P_n\ra \din0, \;\; a\mapsto h_j(a)H_i-h_i(a)H_j,
\end{equation}
is a (well-defined) injection: By (\ref{dus4}),
$ \divv (\phi_{ij}(a)) = (h_ih_j-h_jh_i) (a)=0$, and if $\phi_{ij} (a)=0$ then $h_j(a) H_i = h_i(a) H_j$, and so $ a=0$ since the maps $h_i$ and $h_j$ are bijections. For all $\alpha \in \N^n$ and $i\neq j$, let
\begin{equation}\label{phiij}
\th_{ij}^\alpha := \phi_{ij}(x^\alpha )=x^\alpha ( (\alpha_j+1) H_i-(\alpha_i+1)H_j).
\end{equation}
In particular, $\th_{ij}^0 = H_i-H_j$.
Then
\begin{equation}\label{xaaH}
[x^{\alpha - \alpha_ie_i}x_j\der_i, x_i^{\alpha_i+1}\der_j]=\phi_{ji}(x^\alpha ).
\end{equation}
It is obvious that $\mathfrak{div}_1^0=K\der_1$ and $\mathfrak{div}_1^c=K\der_1+KH_1$.
\begin{lemma}\label{a17Mar13}
Let $n\geq 2$. Then
\begin{enumerate}
\item $\din0 = \bigoplus_{i=1}^{n-1} \phi_{i,i+1} (P_n)$.
\item The set of elements $\th_i^\alpha := \phi_{i,i+1}(x^\alpha ) = x^\alpha ((\alpha_{i+1}+1) H_i-(\alpha_i+1)H_{i+1})$, where $i=1, \ldots , n-1$ and $\alpha \in \N^n$, is a $K$-basis for $\din0$.
\item The set of elements $\th^\alpha_i$ in statement 2 and $x^\beta \der_j$, where $x^\beta \in P_n^{\der_j}$ and $j=1, \ldots , n$, is a $K$-basis for $\divn0$.
    \item The set of elements in statement 3 and $H_i$,  where $i$ is any fixed index in the set $\{ 1, \ldots , n\}$, is a $K$-basis for $\divnc$.
\end{enumerate}
\end{lemma}

{\it Proof}. 1. It is obvious that $R:=\sum_{i=1}^{n-1} \phi_{i,i+1} (P_n)\subseteq \divn0$, see (\ref{dus5}). Recall that $\din0= \divn0\cap D_n'$ and $D'_n=\oplus_{\alpha \in \N^n} x^\alpha \CH_n$. By (\ref{dus5}) and the fact that the  $K$-linear maps $h_1, \ldots , h_n$ are invertible,
$$\din0= R+\din0\cap P_nH_n.$$ By (\ref{dus4}), $\din0\cap P_nH_n=0$. Therefore, $\din0 =R$.

2. Statement 2 follows from statement 1.

3. Statement 3 follows from statement 2 and (\ref{dus2}).

4. Statement 4 follows from statement 3 and the fact that $ \divnc = \divn0\oplus KH_i$, $i=1, \ldots , n$. $\Box $

$\noindent $

Let $\th := x_1\cdots x_n\in P_n$. Then $C_n := \bigoplus_{i\in \N} \th^i \CH_n'$ is an abelian Lie subalgebra of $\divn0$ that contains $\CH_n'$. We will see that $C_n$ is a Cartan subalgebra of the Lie algebras $\divn0$  and $\divnc $(Lemma \ref{da11Mar13}.(3,5)).

By Lemma \ref{a17Mar13}.(2,3), for $n\geq 2$,
\begin{equation}\label{divPH}
C_n=\bigoplus_{i=1}^{n-1}\bigoplus_{m\in \N}K\phi_{i,i+1}(\th^m ),
\end{equation}
\begin{equation}\label{divPH1}
\divn0=\bigoplus_{i=1}^n \bigoplus_{\alpha\in \N^{n,i}}Kx^\alpha\der_i \oplus C_n\oplus \bigoplus_{i=1}^{n-1}\bigoplus_{m\in \N}\bigoplus_{\alpha \in \N^n_d\backslash \{ 0\}}K\phi_{i,i+1}(\th^m x^\alpha ),
\end{equation}
where $\N^n_d:=\cup_{i=1}^n\N^{n,i}=\{ (\alpha_1, \ldots , \alpha_n)\in \N^n \, | \, \alpha_i=0$ for some $i\}$.
 We identify the vector space $\CH_n'=\{ \sum_{i=1}^n \l_iH_i\, | \, \sum_{i=1}^n \l_i=0\}$ with its image in $K^n$ under the $K$-linear injection $\CH_n'\ra K^n$, $ \sum_{i=1}^n \l_iH_i\mapsto  (\l_1, \ldots , \l_n)$. So,
 $$\CH_n'=\{ \l \in K^n \, | \, (\l,\1  ) =\sum_{i=1}^n \l_i=0\}$$ where $\1 :=(1,1, \ldots , 1)$ and $(\l , \mu ) := \sum_{i=1}^n \l_i\mu_i$ is the standard inner product on $K^n$. The dual space $\CH_n'^*:=\Hom_K(\CH_n', K)$ can be identified with the factor space
 $$K^n / K\1 =\{ [\mu ] := \mu +K\1 \; | \; \mu\in K^n\},$$
i.e. $[\mu ] ( \sum_{i=1}^n \l_i H_i) = [\mu ] (\l ) := (\mu , \l ) = \sum_{i=1}^n \mu_i\l_i$.
 By (\ref{divPH1}), the $\CH_n'$-module $\divn0$ is a weight module and the summands in (\ref{divPH1}) are the (nonzero) weight vectors under the adjoint action of $\CH_n'$ on $\divn0$,
 $$ \divn0=\bigoplus_{[\mu ] \in W(\divn0 )}\div0mu$$
where $ \div0mu :=\{ \der \in \divn0 \, | \, [H, \der]=[\mu ] (H) \der$ for all $H\in \CH_n'\}$ is the weight subspace of $\divn0$ that corresponds to the weight $[\mu ] $ and $W=W(\divn0 )$ is the set of weights for $\divn0$. To simplify the notation we identify the set $ \N^n_d$ with its isomorphic copy in the factor vector space $K^n / K\1$ via the map $ K^n \ra K^n / K\1 $, $ \l \mapsto K\1$. So,
\begin{equation}\label{WDnb}
W(\divn0 )  =N^n_d.
\end{equation}
\begin{equation}\label{Dnb1}
\div0mu=\begin{cases}
Kx^\alpha\der_i \oplus  \bigoplus_{i=1}^{n-1}\bigoplus_{m\in \N}\phi_{i,i+1}(\th^m x^{\alpha -e_i+\overline{1} })   & \text{if }[\mu ] =[\alpha - e_i], \alpha\in  \N^{n,i},\\
C_n & \text{if } [\mu ] =0,\\
\bigoplus_{i=1}^{n-1}\bigoplus_{m\in \N}\phi_{i,i+1}(\th^m x^\alpha )& \text{otherwise.}
\end{cases}
\end{equation}
$$\dim_K\, \div0mu =\infty \;\; {\rm for \; all}\;\; [\mu ]\in W(\divn0 ).
$$
The Lie algebra $\divn0 = \ivn0 \oplus \din0$ is the direct sum of its weight $\CH_n'$-submodules with
\begin{equation}\label{Wiv}
W(\ivn0 ) = \coprod_{i=1}^n (\N^{n,i}-e_i), \;\; W(\din0 ) = \N^n_d .
\end{equation}
For $H=\sum_{i=1}^n \l_i H_i\in \CH_n$ and $\alpha \in K^n$, let
$$(H, \alpha )=(\alpha ,H)=\sum_{i=1}^n\alpha_i\l_i.$$ Then, for all $\alpha \in \N^n$, $\divv (x^\alpha H) = (\alpha + \overline{1}, H)$. If, in addition, $H\in \CH_n'$, that is $(H, \overline{1})=0$, then, for all $\alpha \in \N^n$, $\divv (x^\alpha H) = (\alpha +\overline{1}, H) = (\alpha , H)$. It follows that $\din0$ is the direct sum of vector spaces
\begin{equation}\label{Wiv1}
\din0 = \bigoplus_{\alpha \in \N^n} \{ Kx^\alpha H\, | \, (H, \alpha +\overline{1}) =0,  H\in \CH_n\} .
\end{equation}
Let
\begin{equation}\label{Wiv2}
\mathfrak{di}_n'^0 := \bigoplus_{\alpha \in \N^n} \{ Kx^\alpha H\, | \, (H, \alpha ) =0,  H\in \CH_n'\} = \bigoplus_{\alpha \in \N^n_d} \{ K[\th ] x^\alpha H\, | \, (H, \alpha ) =0,  H\in \CH_n'\} .
\end{equation}
Clearly, for $n\geq 2$,  $\mathfrak{di}_n'^0\subset \din0$ and the vector space $\mathfrak{di}_n'^0$ is a left $K[\th ]$-module. We will see shortly that $\mathfrak{di}_n'^0$ is a non-Noetherian Lie algebra (Lemma \ref{a3Apr13}). Notice that $\mathfrak{di}_n'^0=C_2=\oplus_{i\geq 0} K\th^i (H_1-H_2) = K[\th ](H_1-H_2)$ where $\th = x_1x_2$.

$\noindent $

{\bf The commutation relations of the weight vectors in $\divn0$}. By (\ref{divPH1}) and (\ref{WDnb}), there are three types of commutation relations of elements from the weight spaces of the Lie algebra $\divn0$, see (\ref{xadxbd}), (\ref{xadxbd1}) and (\ref{xadxbd2}). For all $x^\alpha\der_i\in P_n^{\der_i}\der_i$ and $ x^\beta \der_j\in P_n^{\der_j}\der_j$,
\begin{equation}\label{xadxbd}
[x^\alpha\der_i, x^\beta \der_j]=\begin{cases}
\phi_{ji}(x^{\alpha+\beta -e_i-e_j})& \text{if }\beta_i\neq 0, \alpha_j \neq 0\\
\beta_ix^{\alpha+\beta -e_i}\der_j& \text{if }\beta_i\neq 0, \alpha_j=0 \\
-\alpha_jx^{\alpha+\beta -e_j}\der_i& \text{if }\beta_i=0, \alpha_j \neq 0\\
0& \text{if }\beta_i=0, \alpha_j =0\\
\end{cases}
\end{equation}
where $\phi_{ji}(x^{\alpha+\beta -e_i-e_j})=x^{\alpha+\beta -e_i-e_j}(\beta_iH_j-\alpha_jH_i)$ (since $\alpha_i=0$ and $\beta_j=0$).

For all elements $x^\alpha H, x^{\alpha'}H'\in \din0$,
\begin{equation}\label{xadxbd1}
[x^\alpha  H, x^{\alpha'} H']=x^{\alpha+\alpha'} ((H,\alpha ')H'-(H', \alpha ) H)\in \din0
\end{equation}
since $((H, \alpha') H'-(H', \alpha ) H, \alpha +\alpha'+\overline{1})= (H, \alpha') ( H',\alpha'+\overline{1})-(H', \alpha )(H, \alpha +\overline{1})=0$.

 For all $x^\beta\in P_n^{\der_i}\der_i$ and $x^\alpha H\in \din0$,
\begin{equation}\label{xadxbd2}
[ x^\beta \der_i, x^\alpha H]=\alpha_i x^{\alpha+\beta-e_i}H -(H, \beta - e_i)x^{\alpha +\beta}  \der_i.
\end{equation}
If, in addition $\alpha_i\neq 0$, then the equality (\ref{xadxbd2}) takes the form
\begin{equation}\label{xadxbd3}
[ x^\beta \der_i, x^\alpha H]=x^{\alpha+\beta-e_i}(\alpha_i H -(H, \beta - e_i)H_i)\in \din0
\end{equation}
since $(\alpha_iH-(H, \beta-e_i)H_i, \alpha +\beta -e_i+\overline{1}) = -(H, \beta -e_i) (H_i, \beta - e_i +\overline{1})=  -(H, \beta -e_i)\cdot 0=0$. By (\ref{Wiv1}) and (\ref{xadxbd1}), $\din0$ is a {\em Lie subalgebra} of $\divn0$ which is not an ideal, by(\ref{xadxbd2}). By (\ref{xadxbd}), $\ivn0$ is not a Lie algebra for $n\geq 2$.

$\noindent $

{\bf The Lie algebra $\mathfrak{di}_n'^0$ is not Noetherian, $n\geq 2$}.  Let $n\geq 2$. The Lie algebra $\mathfrak{di}_n'^0$ is a $K[\th ]$-module, $K[\th ] \mathfrak{di}_n'^0\subseteq \mathfrak{di}_n'^0$, and the Lie bracket on $\mathfrak{di}_n'^0$ is a $k[\th]$-bilinear: for all $p,p'\in K[\th ]$ and $x^\alpha H, x^{\alpha'}H'\in \mathfrak{di}_n'^0$,
\begin{equation}\label{ppHa}
[ px^\alpha H , p'x^{\alpha'} H']=pp'[ x^\alpha H , x^{\alpha'} H'].
\end{equation}
So, the Lie algebra $\mathfrak{di}_n'^0$ is not simple:  $\th^i\mathfrak{di}_n'^0$, $i\in \N$ are distinct ideals of $\mathfrak{di}_n'^0$.
\begin{lemma}\label{a3Apr13}
 For $n\geq2$, the Lie algebras $\din0$ and  $\mathfrak{di}_n'^0$ are not Noetherian.
\end{lemma}

{\it Proof}. If $n=2$ then $\mathfrak{di}_2'^0=C_2$ is an infinite dimensional abelian Lie algebra, hence, non-Noetherian. Let $n\geq 3$. Let $I$ be an ideal of the additive monoid $\N^n$ ($I+\N^n\subseteq I$). Then $\ga' (I):= \bigoplus_{\alpha \in I} \{ x^\alpha H\, | \, $ $(\alpha , H)=0, H\in \CH_n'\}$ is an ideal of the Lie algebra $\mathfrak{di}_n'^0$, by (\ref{xadxbd1}). The map $I\mapsto \ga' (I)$ from the set $\CI (\N^n)$ of ideals of $\N^n$ to the set $\CI (\mathfrak{di}_n'^0)$ of ideals of the Lie algebra $\mathfrak{di}_n'^0$ is an inclusion preserving  injection ($I_1\varsubsetneqq I_2$ implies $\ga' (I_1)\varsubsetneqq\ga' (I_2)$. The set $\CI (\N^n )$ is not Noetherian (with respect to $\subseteq $) hence $\mathfrak{di}_n'^0$ is not a Noetherian Lie algebra.

Similarly, by (\ref{xadxbd1}),  for $n\geq 2$, the map $I\mapsto \ga (I) :=  \bigoplus_{\alpha \in I} \{ x^\alpha H\, | \, $ $(\alpha +\overline{1}, H)=0, H\in \CH_n\}$ from the set $\CI (\N^n)$ to the set $\CI (\din0 )$ of ideals of the Lie algebra $\din0$ is an inclusion preserving injection. Therefore, the Lie algebra $\din0$ is non-Noetherian. $\Box $

$\noindent $

Let $\CG$ be a Lie algebra and $\CH$ be its Lie subalgebra. The {\em centralizer} $C_\CG (\CH ) := \{ x\in \CG \, | \, [ x, \CH ] =0\}$ of $\CH$ in $\CG$ is a Lie subalgebra of $\CG$. In particular, $Z(\CG ) := C_{\CG }(\CG ) $ is the {\em centre} of the Lie algebra $\CG$. The {\em normalizer} $N_\CG (\CH ) :=\{ x\in \CG \, | \, [ x, \CH ] \subseteq \CH\}$ of $\CH$ in $\CG$ is a Lie subalgebra of $\CG$, it is the largest Lie subalgebra of $\CG$ that contains $\CH $ as an ideal.

Let $V$ be a vector space over $K$. A $K$-linear map $\d : V\ra V$
is called a {\em locally nilpotent map} if $V=\bigcup_{i\geq 1} \ker
(\d^i)$ or, equivalently, for every $v\in V$, $\d^i (v) =0$ for
all $i\gg 1$. When  $\d$ is a locally nilpotent map in $V$ we
also say that $\d$ {\em acts locally nilpotently} on $V$. Every {\em nilpotent} linear map  $\d$, that is $\d^n=0$ for some $n\geq 1$, is a locally nilpotent map but not vice versa, in general.   Let
$\CG$ be a Lie algebra. Each element $a\in \CG$ determines the
derivation  of the Lie algebra $\CG$ by the rule $\ad (a) : \CG
\ra \CG$, $b\mapsto [a,b]$, which is called the {\em inner
derivation} associated with $a$. The set $\Inn (\CG )$ of all the
inner derivations of the Lie algebra $\CG$ is a Lie subalgebra of
the Lie algebra $(\End_K(\CG ), [\cdot , \cdot ])$ where $[f,g]:=
fg-gf$. There is the short exact sequence of Lie algebras
$$ 0\ra Z(\CG ) \ra \CG\stackrel{\ad}{\ra} \Inn (\CG )\ra 0,$$
that is $\Inn (\CG ) \simeq \CG / Z(\CG )$ where $Z(\CG )$ is the {\em centre} of the Lie algebra $\CG$ and $\ad ([a,b]) = [
\ad (a) , \ad (b)]$ for all elements $a, b \in \CG$. An element $a\in \CG$ is called a {\em locally nilpotent element} (respectively, a {\em nilpotent element}) if so is the inner derivation $\ad (a)$ of the Lie algebra $\CG$.

$\noindent $

{\bf The Cartan subalgebra $C_n$ of $\divn0$}. A nilpotent Lie subalgebra $C$ of a Lie algebra  $\CG$ such that $ C=N_\CG (C)$  is called a {\em Cartan subalgebra} of $\CG$. We use often the following obvious observation: {\em An  abelian Lie subalgebra that coincides with its centralizer is a maximal abelian Lie subalgebra}.

{\em Example}. $\CH_n$ is a Cartan subalgebra of $D_n$ and $\CH_n = C_{D_n}(\CH_n)$ is a maximal abelian Lie subalgebra of $D_n$.

\begin{lemma}\label{da11Mar13}
Let $n\geq 2$. Recall that   $C_n=\bigoplus_{i\in \N} \th^i \CH_n'$  and $\th = x_1x_2\cdots x_n$. Then
\begin{enumerate}
\item $\CH_n'=\CH_n\cap \divn0$.
\item  $C_n=C_{\divn0}(\CH_n')$  and $ C_n = C_{\divn0} (C_n)$ is a maximal abelian Lie subalgebra of $\divn0$.
\item $C_n$  is a Cartan subalgebra of $\divn0$.
\item $C_n+\CH_n= C_{\divnc} (\CH_n')$ and $C_n = C_{\divnc} (C_n)$ is a maximal abelian Lie subalgebra of $\divnc$.
    \item $N_{\divnc}(C_n) = C_n +\CH_n= N_{\divnc}(C_n + \CH_n)$ is a solvable but not nilpotent Lie algebra, and so $C_n +\CH_n$  and $C_n$ are not a Cartan subalgebra of $\divnc$.
    \item $\CH_n$ is a Cartan subalgebra of $\divnc$ and $\CH_n =  C_{\divnc} (\CH_n)$ is a maximal abelian Lie subalgebra of $\divnc$.
    \end{enumerate}
\end{lemma}

{\it Proof}. 1. Let $H=\sum_{i=1}^n \l_i H_i\in \CH_n$. Then $H\in \CH_n\cap \divn0$ iff $\sum_{i=1}^n \l_i=0$ iff $H\in \CH_n'$.

2. The fact that $C_n =C_{\divn0} (\CH'_n) $ follows from (\ref{Dnb1}). By (\ref{xadxbd1}), $C_n$ is an abelian Lie subalgebra of $\divn0$. Then the inclusion $\CH_n'\subseteq C_n = C_{\divn0} (\CH_n')$ implies the inclusions
$$ C_n = C_{\divn0} (\CH_n')\supseteq C_{\divn0} (C_n) \supseteq C_n.$$
Therefore, $C_n=C_{\divn0} (C_n)$ is a maximal abelian Lie subalgebra of $\divn0$.

3. Let $N:=N_{\divn0}(C_n)$. By statement 2, $C_n\subseteq N$. The $\CH_n'$-module $\divn0$ is weight and $ C_n =  C_{\divn0}(\CH_n')$ is the zero weight component of $\divn0$. The inclusion $ [N, \CH_n']\subseteq C_n$ implies $N\subseteq C_n$, and so $N = C_n$.

4. Notice that $\divnc = \divn0 \oplus KH_1$, $H_1\in C_{\divnc}(\CH_n')$ and $\CH_n' \subseteq \divn0$. By statement 2,
$$ C_{\divnc}(\CH_n')  =  C_{\divn0}(\CH_n')\oplus KH_1= C_n \oplus KH_1=C_n+\CH_n.$$ Now, the inclusion $\CH_n'\subseteq C_n$ implies the inclusions
$$C_n \oplus KH_1=C_{\divnc}(\CH_n') \supseteq C_{\divnc}(C_n) \supseteq C_n.$$ Hence $C_{\divnc}(C_n) =C_n\oplus C_{\divnc}(C_n) \cap KH_1= C_n$ and so $C_n$ is a maximal abelian Lie subalgebra of $\divnc$.

5. Let $N^c:= N_{\divnc}(C_n)$. By statement 4, $C_n\subseteq N^c$. The $\CH_n'$-module $\divnc = \divn0\oplus KH_1$ is weight and $C_n\oplus KH_1= C_{\divnc}(\CH_n')$ is the zero weight component of $\divnc$. The inclusion $[N^c, \CH_n']\subseteq C_n$ implies the inclusion $[N^c, \CH_n']\subseteq C_n\oplus KH_1= C_{\divnc}(\CH_n')$, and so $N^c\subseteq C_n\oplus KH_1$. Since $C_n\subseteq N^c$, it follows that $N^c=C_n \oplus N^c\cap KH_1= C_n+KH_1$, i.e. $C_n$ is not a Cartan subalgebra of $\divnc$.

Clearly, $C_n + \CH_n\subseteq \CN := N_{\divnc}(C_n+\CH_n)$. Since $ \CH_n'\subseteq C_n+\CH_n$, $[\CH_n', \CN ] \subseteq C_n+\CH_n = C_{\divnc}(\CH_n')$, hence $\CN \subseteq C_n + \CH_n$ (as $C_n + \CH_n$ is the zero weight  component of the weight $\CH_n'$-module $\divnc$), i.e. $\CN = C_n + \CH_n$. The Lie algebra $C_n +\CH_n$ is  solvable but not nilpotent, and so $C_n +\CH_n$ is not a Cartan subalgebra of $\divnc$.

6. Statement 6 follows from the Example above.  $\Box $

$\noindent $



The Lie algebra $\divn0$ is a weight $\CH_n'$-module with respect to the adjoint action. In particular, the action of $\CH_n'$ on $\divn0$ is {\em locally finite dimensional}, i.e. for any $\der \in \divn0$, $ \dim_K(\sum_{i\in \N} \ad (\CH_n' )^i(\der ))<\infty$. We can easily verify that the action of the Cartan subalgebra $C_n$ of $\divn0$ on $\divn0$ is not locally finite dimensional, see (\ref{xadxbd1})  and (\ref{xadxbd2}).

$\noindent $

{\bf $P_n$ is a $D_n$-module}. The polynomial algebra $P_n$ is a (left) $D_n$-module: $D_n \times P_n\ra P_n$, $(\der, p)\mapsto \der *p$. In more detail, if $\der = \sum_{i=1}^n a_i\der_i$ where $a_i\in P_n$ then
$$\der * p = \sum_{i=1}^n a_i\frac{\der p}{\der x_i}.$$
The field $K$ is a $D_n$-submodule of  $P_n$ and
\begin{equation}\label{IkerdiK}
P_n^{\CD_n}:=\bigcap_{i=1}^n \ker_{P_n}(\der_i)= K.
\end{equation}

\begin{lemma}\label{xa11Mar13}
\cite{Bav-Aut-Der-Pol} The $D_n$-module $P_n/K$ is simple with $\End_{D_n}(P_n/K)=K\id  $ where $\id$ is the identity map.
\end{lemma}


{\bf  The $G_n$-module $D_n$}. The Lie algebra $D_n$ is a $G_n$-module,
$$ G_n\times D_n\ra D_n, \;\; (\s , \der ) \mapsto \s (\der ) := \s \der \s^{-1}.$$
Every automorphism $\s \in G_n$ is uniquely determined by the elements
$$x_1':=\s (x_1), \; \ldots , \; x_n':=\s (x_n).$$
Let $M_n(P_n)$ be the algebra of $n\times n$ matrices over  $P_n$. The matrix  $J(\s) := (J(\s )_{ij}) \in M_n(P_n)$, where $J(\s )_{ij} =\frac{\der x_j'}{\der x_i}$,   is called the {\em Jacobian matrix} of the automorphism (endomorphism)  $\s$ and its determinant $\CJ (\s ) :=\det \, J(\s)$ is called the {\em Jacobian} of $\s$. So, the $j$'th column of $J(\s )$ is the {\em gradient} $\grad \, x_j':=(\frac{\der x_j'}{\der x_1}, \ldots , \frac{\der x_j'}{\der x_n})^T$  of the polynomial $x_j'$  where $T$ stands for the transposition. Then the derivations
$$\der_1':= \s \der_1\s^{-1}, \; \ldots , \; \der_n':= \s\der_n\s^{-1}$$ are the partial derivatives of $P_n$ with respect to the variables $x_1', \ldots , x_n'$,
\begin{equation}\label{ddp=dxi}
\der_1'=\frac{\der}{\der x_1'}, \; \ldots , \; \der_n'=\frac{\der}{\der x_n'}.
\end{equation}
Every derivation $\der \in D_n$ is a unique sum $\der = \sum_{i=1}^n a_i\der_i$ where $a_i = \der *x_i\in P_n$. Let  $\der := (\der_1, \ldots , \der_n)^T$ and $ \der' := (\der_1', \ldots , \der_n')^T$. Then
\begin{equation}\label{dp=Jnd}
\der'=J(\s )^{-1}\der , \;\; {\rm i.e.}\;\; \der_i'=\sum_{j=1}^n (J(\s )^{-1})_{ij} \der_j\;\; {\rm for }\;\; i=1, \ldots , n.
\end{equation}
In more detail, if $\der'=A\der $ where $A= (a_{ij})\in M_n(P_n)$, i.e. $\der_i'=\sum_{j=1}^n a_{ij}\der_j$. Then for all $i,j=1, \ldots , n$,
$$\d_{ij}= \der_i'*x_j'=\sum_{k=1}^na_{ik}\frac{\der x_j'}{\der x_k}$$
where $\d_{ij}$ is the Kronecker delta function. The equalities above can be written in the matrix form as  $AJ(\s) = 1$ where $1$ is the identity matrix. Therefore, $A= J(\s )^{-1}$.

For all $\s, \tau \in G_n$,
\begin{equation}\label{Jst=JsJ}
J(\s \tau ) = J(\s )\cdot  \s (J(\tau )).
\end{equation}

By taking  the determinants of both sides of (\ref{Jst=JsJ}), we have a similar equality of the Jacobians: for all $\s , \tau \in G_n$.
\begin{equation}\label{Jst=JsJ1}
\CJ (\s \tau ) = \CJ(\s )\cdot  \s (\CJ(\tau )).
\end{equation}

{\bf Properties of the divergence}.  Recall some of the properties of the divergence map
$$ \divv : D_n\ra P_n , \;\; \der =\sum_{i=1}^n a_i\der_i \mapsto \divv (\der ) = \sum_{i=1}^n \frac{\der a_i}{\der x_i}.$$

(div-i) $\divv$ is a $K$-linear map which is a surjection.

$\noindent $

(div-ii) For all $a\in P_n$ and $\der \in D_n$, $\divv (a\der ) = a\divv (\der ) + \der (a)$.

$\noindent $

(div-iii) For all $\der, \d \in D_n$, $\divv ([\der , \d ]) = \der (\divv (\d ))- \d (\divv (\der ))$.

$\noindent $

(div-iv) Let $a_1, \ldots , a_n\in P_n$; $\s :P_n\ra P_n$, $x_i\mapsto a_i$, $i=1, \ldots , n$ and $ \CJ (a_1, \ldots , a_n):=\CJ (\s )$ be the Jacobian of $\s $. Then (Proposition 2.3.2, \cite{Nowicki-book})
$$ \der * \CJ (a_1, \ldots , a_n)= -\CJ (a_1, \ldots , a_n)\divv (\der ) + \sum_{i=1}^n \CJ (a_1, \ldots , \der *a_i, \ldots ,  a_n).$$

$\noindent $

(div-v) (Theorem 2.5.5, \cite{Nowicki-book}) If $D_n=\bigoplus_{i=1}^n P_n\der_i'$ and $\der_1', \ldots , \der_n'$ commute then
$$ \divv (\der_1')=\cdots = \divv (\der_n')=0.$$

$\noindent $

(div-vi)  Let $\s :P_n\ra P_n$, $x_i\mapsto x_i'$, $i=1, \ldots , n$, be an automorphism of $P_n$. Then $\der_1':= \s \der_1\s^{-1}, \ldots , \der_n':= \s \der_n\s^{-1}$ are commuting derivations of $P_n$ such that $D_n = \bigoplus_{i=1}^n P_n\der_i'$ (by (\ref{dp=Jnd}) and since $J(\s)^{-1} \in \GL_n (P_n)$). By (\ref{dp=Jnd}), $\der_i = \sum_{i=1}^n (J(\s )^{-1})_{ij}\der_j$. By (div-v),
$$ \sum_{j=1}^n \der_j*(J(\s )^{-1})_{ij}=0 \;\; {\rm for } \;\; i=1, \ldots , n.$$


{\bf The divergence commutes with polynomial automorphisms}. The following known  theorem shows that the divergence commutes with automorphisms, i.e. the divergence map $\divv : D_n \ra P_n$ is a $G_n$-module homomorphism. We give a short proof.

\begin{theorem}\label{B16Mar13}
For all $\s \in G_n$ and $\der \in D_n$,
$$ \divv (\s (\der ))= \s (\divv (\der )).$$
\end{theorem}

{\it Proof}. Let $\der = \sum_{i=1}^n a_i\der_i$ where $a_i\in P_n$. Then
$ \der' = \s \der \s^{-1} =\sum_{i=1}^n \s (a_i)  \der_i'$ where, by (\ref{dp=Jnd}), $\der_i' = \sum_{j=1}^n (J(\s )^{-1})_{ij}\der_j$. Now, by (div-vi),
\begin{eqnarray*}
 \divv (\der') &=& \sum_{i,j=1}^n \der_j*((J(\s)^{-1})_{ij}\s (a_i)) = \sum_{i=1}^n (\sum_{j=1}^n \der_j *(J(\s )^{-1})_{ij})\cdot \s (a_i))+ \sum_{i=1}^n \sum_{j=1}^n (J(\s)^{-1})_{ij} \der_j*\s (a_i) \\
 &=&  \sum_{i=1}^n \der_i'*\s (a_i) = \sum_{i=1}^n \s \der_i \s^{-1}\s (a_i) = \s ( \sum_{i=1}^n \der_i (a_i)) = \s (\divv (\der )). \;\; \Box
\end{eqnarray*}

\begin{theorem}\label{11Mar13}
$\mG_n = G_n$.
 \end{theorem}

The above theorem was announced in \cite{Rudakov-1986} where a sketch of a proof is given, it can also be deduced from \cite{Skryabin-RegLieDer-1988}.
 A proof of the above theorem is given in \cite{Siebert-1996} and in (\cite{Osborn-AutWitt-1997}, Theorem 3.6),  a different approach and a short proof is given in
\cite{Bav-Aut-Der-Pol}. Some generalizations are given in \cite{Grabowski-1978}, \cite{Dokovic-Zhao-GenCarW-1997} and \cite{Grabowski-Poncin-2004}.

\begin{corollary}\label{a20Mar13}
For all $\s \in \mG_n$ and  $\der \in D_n$,
$ \divv (\s (\der )) = \s (\divv (\der )).$
\end{corollary}
{\it Proof}. The statement follows from Theorem \ref{11Mar13} and Theorem \ref{B16Mar13}. $\Box $

$\noindent $

By (div-iii), $\divn0$ and $\divnc$ are Lie subalgebras of $D_n$, $\divn0$ is an ideal of $\divnc$.

\begin{corollary}\label{aB16Mar13}
The Lie subalgebras $\divn0$ and $\divnc$ of $D_n$ are also $G_n$-submodules of $D_n$.
\end{corollary}

{\it Proof}. This follows from Theorem \ref{B16Mar13}. $\Box $

$\noindent $

By (div-iii),
\begin{equation}\label{divDnPn}
0\ra \divn0\ra D_n \stackrel{\divv}{\ra}P_n\ra 0
\end{equation}
is the short exact sequence of (left) $\divn0$-modules and $\divnc$-modules, i.e.
\begin{equation}\label{divcom}
\divv ([\der , \d ])= \der*\divv (\d )
\end{equation}
for all    $\der\in \divnc$ and $\d \in D_n$. So,
\begin{equation}\label{divDnPn1}
0\ra \divn0\ra \divnc \stackrel{\divv}{\ra}K\ra 0
\end{equation}
is the  short exact sequence of (left) $\divn0$-modules/$\divnc$-modules, and so
\begin{equation}\label{divDnPn2}
 \divnc =\divn0\oplus KH_i\;\; {\rm for}\;\; i=1, \ldots , n,
\end{equation}
 since $\divv (H_i) =1$ for all $i=1, \ldots , n$.

$\noindent $

{\bf The maximal abelian Lie subalgebra $\CD_n$ of $\divn0$ and $\divnc$}.


\begin{lemma}\label{db11Mar13}
\begin{enumerate}
\item $C_{\divn0}(\CD_n) = C_{\divnc}(\CD_n) = \CD_n$ and so $\CD_n$ is a maximal abelian Lie subalgebra of $\divn0$ and $\divnc$.
\item (\cite{Bav-Aut-Der-Pol})  $\Fix_{G_n} ( \CD_n) = \Sh_n$.
\item For $n\geq 2$, $\divn0$ is a faithful $G_n$-module, i.e. the group homomorphism $G_n\ra \mfGn$, $\s \mapsto \s : \der\mapsto \s \der \s^{-1}$, is a monomorphism. For $n=1$, the group $\Sh_1$ is the kernel of the group homomorphism $G_1\ra \mathbf{G}_1$.
  \item $\divnc$ is a faithful $G_n$-module, i.e. the group homomorphism $G_n\ra \mfGnc$, $\s \mapsto \s : \der\mapsto \s \der \s^{-1}$, is a monomorphism.
      \item $\Fix_{G_n} (\CD_n\oplus \CH_n') = \{ e\}$ for $n\geq 2$.
    \end{enumerate}
\end{lemma}

{\it Proof}. 1. Since $C_{D_n}(\CD_n) = \CD_n$, \cite{Bav-Aut-Der-Pol}, we must have $C_{\divn0} (\CD_n) = C_{\divnc} (\CD_n) =\CD_n$, and so $\CD_n$ is a maximal abelian Lie subalgebra of $\divn0$ and $\divnc$.

3 and 5. By Corollary \ref{aB16Mar13}, $\divn0$ is a $G_n$-submodule of $\CD_n$. So, the group homomorphism in statement 3 is well-defined. The case $n=1$ is obvious since $G_1=\{ x\mapsto \l x+\mu \, | \, \l \in K^*, \mu \in K\}$ and $\mathfrak{div}_1^0=K\der$.

So, let $n\geq 2$,  and $\s \in \Fix_{G_n} (\CD_n \oplus \CH_n')$. Then $\s \in \Fix_{G_n}(\CD_n) = \Sh_n$, by statement 2. So, $\s  (x_1)=x_1+\l_1, \ldots , \s (x_n) = x_n+\l_n$ where $\l_i\in K$. Then for all $i\neq j$,
$$ H_i-H_j= \s (H_i-H_j) = (x_i+\l_i) \der_i - (x_j +\l_j) \der_j = H_i- H_j+\l_i \der_i -\l_j\der_j. $$
So, $\l_1=\cdots =\l_n=0$. This means that $\s =e$. So, $\Fix_{G_n} (\CD_n \oplus \CH_n') = \{ e\}$ and $\divn0$ is a faithful $G_n$-module.

4.  By Corollary \ref{aB16Mar13}, $\divnc$ is a $G_n$-submodule of $D_n$. So, the group homomorphism in statement 4 is well-defined. Now, statement 4 follows from statement 3 for $n\geq 2$. For $n=1$, statement 4 is obvious as $\mathfrak{div}_1^c=K\der_1+KH_1$ and $G_1= \{ x\mapsto \l x+\mu \, | \, \l \in K^*, \mu \in K\}$.  $\Box $


\begin{lemma}\label{a16Mar13}
$[ \divnc, \divnc ] = \divn0$.
\end{lemma}

{\it Proof}. The statement is obvious for $n=1$ as $\mathfrak{div}_1^c=K\der_1+KH_1$ and $\mathfrak{div}_1^0=K\der_1$. So, let $n\geq 2$.
 By (div-iii), $\divv ([ \divnc, \divnc ])=0$, and so $[ \divnc, \divnc ]\subseteq \divn0$. Recall that $\divn0 = \ivn0 +\din0$ (see (\ref{dus2})), $\ivn0 = \oplus_{i=1}^nP^{\der_i}\der_i$, $[x^\alpha \der_i , H_i] = x^\alpha \der_i$ for all $x^\alpha \in P_n^{\der_i}$  and $i=1, \ldots , n$. Hence $\ivn0\subseteq [ \divnc, \divnc ]$.

Finally, for all $a\in P_n$ and $i\neq j$,
$$ [H_k, h_j(a)H_i- h_i(a) H_j] = h_j(H_k*a) H_i-h_i (H_k*a) H_j.$$
Hence, $\{ h_j(x^\alpha)H_i- h_i(x^\alpha) H_j\, | \, 0\neq \alpha \in \N^n, i\neq j\} \subseteq [ \divnc, \divnc ]$. Finally, $H_i-H_j= [ x_i\der_j , x_j\der_i]\in [ \divnc, \divnc ]$ for all $i\neq j$. Now, by Lemma \ref{a17Mar13}.(3), $[ \divnc, \divnc ] = \divn0$. $\Box $

The following lemma is well-known and easy to prove.
\begin{lemma}\label{c11Mar13}
\begin{enumerate}
\item $D_n$ is a simple Lie algebra.
\item $Z(D_n)=\{ 0\}$.
\item $[D_n, D_n]=D_n$.
\end{enumerate}
\end{lemma}


The next lemma is also known but we give a short elementary proof.

\begin{lemma}\label{dc11Mar13}
\begin{enumerate}
\item $\divn0$ is a simple Lie algebra.
\item $Z(\divn0 ) =\begin{cases}
\divn0 & \text{if }n=1,\\
0 & \text{if }n\geq 2.
\end{cases}$
\item $[\divn0 , \divn0 ] =\begin{cases}
0& \text{if }n=1,\\
\divn0 & \text{if }n\geq 2.
\end{cases}$
\end{enumerate}
\end{lemma}

{\it Proof}. All three statements are obvious if $n=1$ since $\mathfrak{div}_1^0 = K\der$. So, we assume that $n\geq 2$.

1.  Let $ 0\neq a\in \divn0$ and $ \ga = (a)$ be the ideal of the Lie algebra $\divn0$ generated by the element $a$. We have to show that $\ga = \divn0$. Using the inner derivations $\d_1, \ldots , \d_n$ and $\ad (H_i-H_j)$ (where $i\neq j$)  of the Lie algebra $\divn0$  we see that $\der_i\in \ga$ for some $i$. Then $\CD_n\subseteq \ga$ since $ \der_j = [ \der_i , x_i\der_j] $ for all $i\neq j$. Then $ \ivn0\subseteq \ga$ since $ x^\alpha \der_i =[ \der_j , (\alpha_j+1)^{-1} x^{\alpha+e_j}\der_i ] $ for all $i\ne j$ and $ x^\alpha \in P_n^{\der_i}$. Then $\CH_n'\subseteq \ga$ since $ [ x_i \der_j , x_j \der_i] = H_i-H_j$ for all $i\neq j$.  Using the commutation relations
\begin{equation}\label{Hsthj}
[H_s-H_t, h_j(x^\alpha ) H_i - h_i(x^\alpha ) H_j]= (\alpha_s-\alpha_t) (h_j(x^\alpha ) H_i - h_j(x^\alpha ) H_j)
\end{equation}
we see that  $ \{ h_j(x^\alpha ) H_i - h_j(x^\alpha ) H_j \, | \, i\neq j, \alpha \in \N^n , \alpha \neq  (p,p, \ldots , p), p=1, 2,  \ldots \}\subseteq \ga$. Let $\th = x_1\cdots x_n$. Finally, by Lemma \ref{a17Mar13}.(3),
 $\ga = \divn0$ since for all $p=1, 2, \ldots $ and $i\neq j$,
\begin{equation}\label{Hsthj1}
[x_j\der_i , h_j(x_ix_j^{-1}\th^p ) H_i - h_i( x_ix_j^{-1} \th^p ) H_j] = (p+2) (h_j(\th^p) H_i - h_i(\th^p) H_j).
\end{equation}
In more detail,
\begin{eqnarray*}
 {\rm LHS} &=&  [ x_j\der_i , p x_ix_j^{-1} \th^p H_i - (p+2) x_ix_j^{-1} \th^p H_j] = p(p+2) \th^p H_i - (p+2) (p+1) \th^p H_j\\
 &+&  (p+2) \th^p H_i= (p+2) (p+1) (\th^p H_i -\th^p H_j)\\
 & = &(p+2) (h_j(\th^p) H_i - h_i(\th^p) H_j).
\end{eqnarray*}

2. Statement 2 follows from statement 1.

3. Statement 3 follows from statement 1. $\Box $

$\noindent $


For any $\alpha = (\alpha_1, \ldots , \alpha_n) \in \N^n$ and any $i\neq j$, we have $ \alpha = \alpha_ie_i +\alpha_j e_j+\beta $ where  $\beta \in\N^n$ with $\beta_i = \beta_j=0$. Then
\begin{equation}\label{taijx}
\th_{ij}^\alpha = [ x_i^{\alpha_i+1}\der_j , x^{\alpha - \alpha_ie_i +e_j}\der_i ].
\end{equation}
Indeed, let $c$ be the commutator. Then
$$ c= x^\beta [ x_i^{\alpha_i +1} \der_j , x_j^{\alpha_j+1}\der_i] = x^\beta \th_{ij}^{\alpha_ie_i+ \alpha_j e_j} = \th_{ij}^\alpha .$$ By (\ref{taijx}),
\begin{equation}\label{taijx1}
\din0 \subseteq [ \ivn0 , \ivn0 ].
\end{equation}
The inner derivations $\d_1=\ad (\der_1), \ldots , d_n=\ad (\der_n)$ of the Lie algebra $D_n$ are commuting and locally nilpotent. The Lie algebra $D_n$ is a union of vector spaces
$$ D_n = \bigcup_{i\geq o}D_{n,i}\;\; {\rm where}\;\; D_{n,i} := \{ \der \in D_n \, | \, \d^\alpha (\der ) =0\;\; {\rm for\; all}\;\; \alpha \in \N^n\;{\rm such \; that}\;\; |\alpha | =i+1\}$$
where $\d^\alpha = \prod_{i=1}^n \d_i^{\alpha_i}$. Clearly, $D_{n,i} = \oplus_{j=1}^n P_{n,i}\der_j$ where $P_{n,i}:=\{ p\in P_n\, | \, \deg (p)\leq i\}$ and $[D_{n,i}, D_{n,j}]\subseteq D_{n, i+j-1}$ for all $i,j\geq 0$ where $D_{n,-1}:=0$. The inner derivations $\d_1, \ldots ,\d_n$ are also commuting locally nilpotent derivations of the Lie algebras $\divn0$ and $\divnc$. For each $i\in \N$, let $\mathfrak{div}_{n,i}^0:= \divn0\cap D_{n,i}$ and $\mathfrak{div}_{n,i}^c:= \divnc\cap D_{n,i}$. Then, for all $\alpha\in \N^n$ and $i=1, \ldots , n-1$,
\begin{equation}\label{tiadi}
\th_i^\alpha \in \mathfrak{div}_{n,|\alpha | +1}^0\backslash \mathfrak{div}_{n,|\alpha |}^0.
\end{equation}

 {\bf The automorphisms $s_{i,i+1}$, $i=1, \ldots , n-1$ $(n\geq 2)$}. The automorphism $s= s_{i,i+1}$ of the polynomial algebra $P_n$ that swaps the variables $x_i$ and $x_{i+1}$ leaving the rest of the variables untouched ($s(x_i) = x_{i+1}$ and $s(x_{i+1}) = x_i$) extends uniquely  to an automorphism of the Lie algebra $D_n$. Clearly, $s$ is also an automorphism of the Lie algebras $\divn0$ and $\divnc$ (Theorem \ref{B16Mar13}). In particular,
 $$ s(\der_i) = \der_{i+1}, \; s(\der_{i+1}) = \der_i , \; s(H_i) = H_{i+1}\; {\rm and}\; s(H_{i+1}) = H_i.$$
Therefore, for all $\alpha \in \N^n$,
\begin{equation}\label{siia}
s_{i,i+1}(\th_i^\alpha ) = - \th_i^{s_{i,i+1}(\alpha ) }
\end{equation}
where, for $\alpha = (\alpha_1, \ldots , \alpha_n)$, $s_{i,i+1}(\alpha ) = (\alpha_1, \ldots , \alpha_{i-1}, \alpha_{i+1}, \alpha_i, \alpha_{i+2},\ldots , \alpha_n)$.
For $i=1, \ldots , n-1$; $j\in \N$ and $\alpha \in \N^n$ with $\alpha_i\geq 1$,
\begin{equation}\label{ditai}
[x_{i+1}^j\der_i, \th_i^\alpha]=(\alpha_i+1) \th_i^{\alpha - e_i +je_{i+1}}.
\end{equation}
{\it Proof}. By (\ref{xadxbd3}),
\begin{eqnarray*}
{\rm LHS} &=& [x_{i+1}^j\der_i,x^\alpha ((\alpha_{i+1}+1) H_i-(\alpha_i +1) H_{i+1}] \\
 &=&x^{\alpha -e_i +je_{i+1}}(\alpha_i ((\alpha_{i+1}+1)H_i-(\alpha_i +1) H_{i+1})-((\alpha_{i+1}+1) H_i - (\alpha_i+1) H_{i+1}, je_{i+1}-e_i)H_i)\\
 &=&(\alpha_i+1)x^{\alpha -e_i +je_{i+1}}((\alpha_{i+1}+1+j)H_i-\alpha_iH_{i+1})\\
 &=& (\alpha_i+1) \th_i^{\alpha - e_i +je_{i+1}}.\;\;\; \Box
\end{eqnarray*}

For $i=1, \ldots , n-1$; $j\in \N$ and $\alpha \in \N^n$ with $\alpha_{i+1}\geq 1$,
\begin{equation}\label{ditai1}
[x_i^j\der_{i+1}, \th_i^\alpha]=(\alpha_{i+1}+1) \th_i^{\alpha + je_i -e_{i+1}}.
\end{equation}
{\it Proof}. By applying the automorphism $s=s_{i,i+1}$ to the equality (\ref{ditai}) and using (\ref{siia}) we obtain (\ref{ditai1}):
\begin{eqnarray*}
 [x_i^j\der_{i+1}, \th_i^\alpha]&=& -s([x_{i+1}^j\der_i, \th_i^{s(\alpha )} ]=-s((\alpha_{i+1}+1)\th_i^{s(\alpha )-e_i+je_{i+1}} \\
 &=&(\alpha_{i+1}+1) \th_i^{s(s(\alpha ) -e_i+je_{i+1})}\\
 &=&(\alpha_{i+1}+1) \th_i^{\alpha + je_i -e_{i+1}}. \;\;\; \Box
  \end{eqnarray*}

By (\ref{ditai}) and (\ref{ditai1}), for $i=1, \ldots , n-1$; $j=1, \ldots , n$ and $\alpha \in \N^n$.

\begin{equation}\label{ditai4}
[\der_j, \th_i^\alpha ] = \begin{cases}
(\alpha_j+1) \th_i^{\alpha - e_j}& \text{if }j=i,i+1\; {\rm and} \; \alpha_j\geq 1,\\
(\alpha_{i+1}+1)x^\alpha \der_i& \text{if }j=i, \alpha_i=0,\\
-(\alpha_i+1)x^\alpha \der_{i+1}&  \text{if }j=i+1, \alpha_{i+1}=0,\\
\alpha_j\th_i^{\alpha - e_j}& \text{otherwise}.
\end{cases}
\end{equation}

\begin{proposition}\label{dB11Mar13}
Let $n\geq 2$. Then
\begin{enumerate}
\item $\Fix_{\mfGn} (\CD_n\oplus \CH_n') = \{ e\}$.
\item Let $\s ,\tau \in \mfGn$. Then $\s = \tau$ iff $\s (\der_i)= \tau ( \der_i)$ for $i=1, \ldots , n$ and $\s (H_j-H_{j+1}) = \tau (H_j-H_{j+1})$ for $j=1, \ldots , n-1$.
\item $\Fix_{\mfGn} (\der_1, \ldots , \der_n) = \Sh_n$.
\end{enumerate}
\end{proposition}

{\it Proof}. 1. Let $\s \in F:= \Fix_{\mfGn} (\CD_n\oplus \CH_n')$. We have to show that $\s =e$. Let $\CF := \Fix_{\divn0} (\s ) := \{ \der \in \divn0 \, | \, \s (\der )= \der \}$. We have to show that $\CF = \divn0$.

(i) $\CD_n +\CH_n'\subseteq \CF$: Obvious.

(ii) $\ivn0\subseteq \CF$: We have to show that $P_n^{\der_i}\der_i \subseteq \CF$ for all $i=1, \ldots , n$. Fix $i$ and $x^\alpha \der_i\in P_n^{\der_i}\der_i$. If $\alpha =0$ then $\der_i\in \CD_n \subseteq \CF$, by (i). We use induction on $|\alpha | = \sum_{i=1}^n \alpha_i$ to prove the statement. The case $|\alpha |=0$ is obvious. So, let $|\alpha |\geq 1$.
  The element $x^\alpha \der_i$ has weight $[\alpha - e_i]$. The weight subspace $\mathfrak{div}_{n,[\alpha - e_i]}^0$ is $\s$-invariant ($\s (\mathfrak{div}_{n,[\alpha - e_i]}^0)=\mathfrak{div}_{n,[\alpha - e_i]}^0$) and the set $\{ x^\alpha \der_i, \th_j^{\alpha -e_i+k\overline{1}}\, | \, j=1, \ldots , n-1; \, k=1,2,\ldots \}$ is its $K$-basis, by (\ref{Dnb1}). The element $x^\alpha \der_i$ belongs to the $\s$-invariant subspace $N_{|\alpha | }:= \mathfrak{div}_{n,|\alpha |}^0 (= \divn0 \cap D_{n, |\alpha |})$. By (\ref{tiadi}) and (\ref{ditai4}),
  $$ V:= N_{|\alpha |}\cap \mathfrak{div}_{n,[\alpha - e_i]}^0= Kx^\alpha \der_i$$ is a $\s$-invariant subspace.  Therefore, $\s ( x^\alpha \der_i) = \l_{\alpha , i} x^\alpha \der_i$ for all $x^\alpha \der_i\in P_n^{\der_i}\der_i$ and $i=1, \ldots , n$ for some $\l_{\alpha , i}\in K$. Clearly, $\l_{0,i}=1$ for all $i=1, \ldots , n$ since $\s (\der_i) = \der_i$. Applying the automorphism $\s$ to the relations $[\der_j , x^\alpha \der_i] = \alpha_j x^{\alpha - e_j} \der_i$ yields the equalities $\alpha_j( \l_{\alpha , i} - \l_{\alpha - e_j, i})=0$. If $\alpha_j\neq 0$ then $\l_{\alpha , i} = \l_{\alpha - e_j, i}=1$ (by  induction). Hence, $\ivn0\subseteq \CF$.

(iii) $\din0 \subseteq \CF$: By (ii) and (\ref{taijx1}).

Therefore, $\divn0 = \din0 +\ivn0 \subseteq \CF$, and so $\divn0 = \CF$, as required.

2. Statement 2 follows from statement 1.

3. Clearly, $\Sh_n \subseteq F: = \Fix_{\mfGn}(\der_1, \ldots , \der_n)$. Let $\s \in F$ and $H_{i,j}':=\s (H_{i,j})$ where $H_{i,j}:= H_i-H_j$ for $i\neq j$. Then
$$ [ \der_k , H_{i,j}'-H_{i,j}]= \s ([ \der_k, H_{i,j}]) - [\der_k , H_{i,j}] = [\der_k, H_{i,j}] - [\der_k , H_{i,j}]=0$$
since $[\CD_n , \CH_n']\subseteq \CD_n$. Then $d_{ij} := H_{i,j}'-H_{i,j}\in C_{\divn0} (\CD_n) = \CD_n$, and so $d_{ij}=\sum_{k=1}^n \l_{ij}^k \der_k$ for  some $\l_{ij}^k\in K$.

If $n=2$ then $H_{1,2}'=H_1-H_2+\l_1\der_1-\l_2\der_2= (x_1+\l_1) \der_1-(x_2+\l_2)\der_2 = s_\l (H_1-H_2)$ where $s_\l \in \Sh_2$ and $\l = (\l_1, \l_2 )$. Hence, $s_\l^{-1}\s   (H_1-H_2) = H_1-H_2$, i.e.
 $s_\l^{-1}\s \in \Fix_{\mathbf{G}_2}(\CD_2+\CH_2') = \{ e\}$, and so $\s =s_\l \in \Sh_2$.

 Suppose that $n>2$. If $n=3$ then up to  action of $\Sh_3$, we may assume that $ H_{1,2}'= H_{1,2} +\l \der_3$ and $ H_{2,3}'= H_{2,3} +\mu \der_1+\nu \der_2$ for some $\l , \mu , \nu \in K$. Then $0=[H_{1,2}', H_{2,3}'] = -\mu \der_1+\nu \der_2+\l \der_3$ and so $\l = \mu = \nu =0$.

 If $n\geq 4$ then for any four distinct numbers $i,j,k,l\in \{ 1,2, \ldots , n\}$ we have the equality
 $$ 0= [ H_{i,j}', H_{k,l}'] = [ H_{i,j}, d_{kl}]-[H_{k,l},  d_{ij}] = -\l_{kl}^i\der_i +\l_{kl}^j\der_j +\l_{ij}^k\der_k - \l_{ij}^l\der_l.$$
Therefore, $\l^i_{kl}=0$ for all distinct $i$, $j$ and $k$. Then using $\Sh_n$, we may assume that
$$H_{1,2}'= H_{1,2}, \;\; H_{2,3}'= H_{2,3}+\l_{23}^2\der_2, \;\;
 H_{3,4}'= H_{3,4}+\l_{34}^3\der_3, \ldots  ,  H_{n-1,n}'= H_{n-1,n}+\l_{n-1,n}^{n-1}\der_{n-1}.$$
Then, for $i=2, \ldots , n-1$, $0=[H_{i-1, i}', H_{i,i+1}']=\l_{i, i+1}^i\der_i$, i.e. $\l_{i, i+1}^i=0$ and $H_{i,i+1}'=H_{i,i+1}$. This means that $ s_\l \s \in \Fix_{\mfGn} (\CD_n+\CH_n')=\{ e\}$ (statement 2) for some $s_\l \in \Sh_n$, and so $\s = s_\l^{-1}\in \Sh_n$. $\Box $


\begin{lemma}\label{dc13Mar13}
Let $\s \in \mfGn$,  $\der_1':=\s (\der_1), \ldots , \der_n':= \s (\der_n)$;  $\d_1':=\ad (\der_1'), \ldots , \d_n':= \ad (\der_n')$ and $\CD'$ be  the subalgebra of $\End_K(P_n)$ generated by the linear  maps $\der_1',\ldots , \der_n'$. Then
\begin{enumerate}
\item $\der_1', \ldots , \der_n'$ are commuting, locally nilpotent derivations of $P_n$.
\item $\bigcap_{i=1}^n\ker_{P_n}(\der_i')=K$.
\item {\rm (Embedding trick)} For $n\geq 2$, the short sequence of $\CD'$-modules $0\ra K\ra P_n\stackrel{\D'}{\ra} {\divn0}^{n\choose 2}$ is exact where $\D' (p) :=\prod_{i<j} (\der_i'(p)\der_j'-\der_j'(p)\der_i')$. In particular, for all $p\in P_n$ and $\alpha \in \N^n$,
 $$ \D'\der'^\alpha (p) = \d'^\alpha \D' (p)$$
 where $\der'^\alpha := \prod_{i=1}^n \der_i'^{\alpha_i}$ and   $\d'^\alpha := \prod_{i=1}^n \d_i'^{\alpha_i}$.
 \item Let $\D': P_n\ra D_n^{n\choose 2}$, $p\mapsto \prod_{i<j} (\der_i'(p)\der_j'-\der_j'(p)\der_i')$. Then $\D'\in \Der_{\CD'} (P_n, D_n^{n\choose 2})$.
\end{enumerate}
\end{lemma}

{\it Proof}. 2.   Let $\l \in \bigcap_{i=1}^n \ker_{P_n}(\der_i')$. Then $\divv (\l \der_1') = \l \divv (\der_1')+\der_1'*\l=0$, i.e. $\l \der_1'\in \divn0$ and
 $$\l \der_1' \in C_{\divn0}(\der_1', \ldots , \der_n')=\s (C_{\divn0}(\der_1, \ldots , \der_n))= \s (C_{\divn0}(\CD_n))= \s (\CD_n) = \s (\bigoplus_{i=1}^n K\der_i) = \bigoplus_{i=1}^n K\der_i', $$
since $C_{\divn0}(\CD_n)=\CD_n$, Lemma \ref{db11Mar13}.(1). Then $\l \in K$ since otherwise the infinite dimensional space $\bigoplus_{i\geq 0} K\l^i \der_1'$ would be a subspace of the  finite dimensional space $\s (\CD_n)$.

3. The derivations $\der_1', \ldots , \der_n'$  commute  since $\der_1, \ldots , \der_n$  do. So, $\CD'$ is a commutative algebra. For all $p\in P_n$ and $\der\in \divn0$,
$$ \divv (a\der ) = a\divv (\der ) +\der (a)= \der (a).$$
So, the map $\D'$ is well-defined: $\divv (\der_i'(p)\der_j'-\der_j'(p)\der_i') = (\der_j'\der_i'-\der_i'\der_j')(a)=0$. Recall that the vector spaces $P_n$ and $\divn0$ are left $\divn0$-modules hence they are also left $\CD'$-modules since $\der_1', \ldots , \der_n'\in \divn0$. The map $\D'$ is a $\CD'$-homomorphism since, for all $p\in P_n$,
$$\D'\der_i'(p) = [ \der_i' , \D'(p)], \;\; i=1, \ldots , n.$$
It remains to show that $\ker (\D')= K$. The inclusion $K\subseteq \ker (\D')$ is obvious.

 (i) $\ker (\D') = \{ p\in P_n\, | \, \der_i'(p) \der_j'= \der_j'(p) \der_i'$ {\em for all} $i\neq j\}$: This is obvious. Let $P_n^{\der_i'}:= \ker_{P_n}(\der_i')$.

 (ii) $\ker (\D') \cap  P_n^{\der_i'}= K$ {\em for} $i=1, \ldots , n$: Given $p\in P_n$. By (i), $p\in\ker (\D') \cap  P_n^{\der_i'}$ iff $p\in\ker (\D') \cap  P_n^{\der_j'}$ for $j=1, \ldots , n$ ($
  0=\der_i'(p) \der_j'= \der_j'(p) \der_i' \; \Rightarrow \;  \der_j'(p)=0$) iff
  $$p\in \ker (\D') \cap  P_n^{\der_1', \ldots , \der_n'}=\ker (\D') \cap K = K,$$ by statement 2.
  Suppose that $K':= \ker (\D') \backslash K \neq \emptyset$, we seek a contradiction.

 (iii) {\em If $p\in K'$ then $\der_1'(p) \neq0, \ldots , \der_n'(p)\neq 0$}: This follows from (ii) ($K'=\ker (\D')\backslash K = \ker (\D') \backslash \ker (\D')\cap  P_n^{\der_i'}= \ker (\D') \backslash  P_n^{\der_i'}$).

 (iv) $P_n^{\der_1'}=\cdots =P_n^{\der_n'}$: Fix $p\in K'$. By (iii), $\der_1'(p) \neq0, \ldots , \der_n'(p)\neq 0$. By (i),
 $$ \der_i'(p) \der_j'= \der_j'(p) \der_i'\;\; {\rm for \; all} \;\; i\neq j,$$
 and so $P_n^{\der_1'}=\cdots =P_n^{\der_n'}$.

 (v) $P_n^{\der_1'}=\cdots =P_n^{\der_n'}=K$: This statement follows from (iv) and statement 2.

Suppose that $\ker (\D')\neq K$, we seek a contradiction. Fix an element $p\in \ker (\D') \backslash K$. By statement 2, $\der_j' (p)\neq 0$ for some $j$. For all $i\neq j$,  $\der_i'(p) \der_j'= \der_j'(p) \der_i'$ (since $p\in \ker (\D')$). Hence, $\der_i'(p)\neq 0$ for all $i$. Therefore, for $i=1, \ldots , n$, $\der_i'=f_i\der $ for some $f_i\in P_n$ and a derivation $\der = \sum_{s=1}^n p_s\der_s\in D_n$ with $\gcd (p_1,\ldots , p_n)=1$. For all elements $c_i\in C_i':=C_{\divn0}(\der_i')$, $0=[c_i,\der_i']=c_i(f_i)\der +f_i[c_i,\der ]$, and so $[c_i,\der ]\in P_n\der$, by the choice of $\der$. By Theorem \ref{6Oct13}, the Lie algebra $\divn0$ is generated by the set $C=\sum_{i=1}^n C_{\divn0}(\der_i)$ and also by the set $\s (C)=\sum_{i=1}^n C_i'$. Hence, $[\divn0 , \der ] \subseteq P_n\der$. In particular, for all $i$, $[\der_i, \der ]=\sum_{j=1}^n \der_i(p_j)\der_j\in P_n\der$. Hence, $\der \in \CD_n$, by degree argument. Then $P_n^{\der}\neq K$, and so $P_n^{\der_i'}=P_n^{\der}\neq K$. This contradicts to $(v)$.

4.  The map $\D'$ can be seen as the map
\begin{equation}\label{DpPb}
\D': P_n\ra D_n^{n\choose 2},\;\;  p\mapsto \prod_{i<j} (\der_i'(p)\der_j'-\der_j'(p)\der_i').
\end{equation}
 Notice that $P_n$ and $D_n^{n\choose 2}$ are left $P_n$-modules.

 $\D'\in \Der_{\CD'} (P_n, D_n^{n\choose 2})$, {\em i.e. the map $\D'$ is a $\CD'$-derivation from the polynomial algebra $P_n$ to the left $P_n$-module  $D_n^{n\choose 2}$, i.e. for all} $p,q\in P_n$, $$\D' (pq) = q\D'(p)+p\D'(q):$$
\begin{eqnarray*}
 \D' (pq) & = &  \prod_{i<j} ((q\der_i'(p)+p\der_i'(q))\der_j'-(q\der_j'(p)+p\der_j'(q))\der_i')\\
 &=& q\prod_{i<j} (\der_i'(p)\der_j'-\der_j'(p)\der_i')+p\prod_{i<j} (\der_i'(q)\der_j'-\der_j'(q)\der_i')\\
 &=& q\D'(p)+p\D'(q).
\end{eqnarray*}

 1. 
  The inner derivations $\d_1, \ldots , \d_n$  of the Lie algebra $\divn0$ are commuting and locally nilpotent.  Hence so are  the inner derivations
 $\d_1, \ldots , \d_n'$,
 by statements 2 and 3.  $\Box $




\begin{theorem}\label{A13Mar13}
\cite{Bav-Aut-Der-Pol} Let $\der_1', \ldots , \der_n'$ be commuting, locally nilpotent derivations of the polynomial algebra $P_n$ such that $\bigcap_{i=1}^n \ker_{P_n}(\der_i')=K$. Then there exist polynomials $x_1', \ldots , x_n'\in P_n$ such that
\begin{equation}\label{con*}
\der_i'*x_j'=\d_{ij}\;\; {\rm for}\;\; i,j=1, \ldots , n.
\end{equation}
Moreover, the algebra homomorphism
$$\s : P_n\ra P_n , \;\; x_1\mapsto x_1', \ldots , x_n\mapsto x_n'$$ is an automorphism such that $\der_i'= \s \der_i \s^{-1} = \frac{\der}{\der x_i'}$ for $i=1, \ldots , n$.
\end{theorem}

\begin{corollary}\label{db13Mar13}
Let $\s \in \mfGn$. Then $\tau \s \in \Fix_{\mfGn} (\der_1, \ldots , \der_n)$ for some $\tau \in G_n$.
\end{corollary}

{\it Proof}. By Lemma \ref{dc13Mar13}, the elements $\der_1':= \s (\der_1), \ldots , \der_n':= \s (\der_n)$ satisfy the assumptions of Theorem \ref{A13Mar13}. By Theorem \ref{A13Mar13},  $\der_1':= \tau^{-1}  (\der_1), \ldots , \der_n':= \tau^{-1}  (\der_n)$ for some $\tau \in G_n$. Therefore, $ \tau \s \in \Fix_{\mfGn} (\der_1, \ldots , \der_n)$. $\Box $

$\noindent $

{\bf Proof of Theorem \ref{16Mar13}}. If $n=1$ then $\mathfrak{div}_1^0 = K\der_1$, $\mathbf{G}_1 = \mT^1 \simeq G_1/ \Sh_1$ since $G_1= \mT^1 \ltimes \Sh_1$.

So, let $n\geq 2$. Let $\s \in \mfGn$. By Corollary \ref{db13Mar13}, $\tau \s \in \Fix_{\mfGn} (\der_1, \ldots , \der_n)=\Sh_n$ (Proposition \ref{dB11Mar13}.(3)). Therefore, $\s \in G_n$, i.e. $\mfGn = G_n$. $\Box$


\begin{lemma}\label{a18Mar13}
 $ \Fix_{\mfGnc}(\divn0 ) = \begin{cases}
\Sh_1 & \text{if }n=1,\\
\{ e\} & \text{if }n\geq 2.\\
\end{cases}$
\end{lemma}

{\it Proof}. If $n=1$, $\s (\der ) = \der$ for some $ \s \in\mathbf{G}_1^c$ then applying $\s $ to the equality $ [ \der_1 , H_1 ] = \der_1$ yields $\s (H_1) = H_1+\l \der_1= (x_1+\l ) \der_1= t_\l (H_1)$ for some $\l \in K$ where $t_\l \in \Sh_1$. Hence $\Fix_{\mathbf{G}_1^c}(\mathfrak{div}_1^0)=\Sh_1$.

So, let $n\geq 2$. Let $\s \in F:= \Fix_{\mfGnc}(\divn0 )$, $H_1':= \s (H_1) , \ldots , H_n':= \s (H_n)$. By (\ref{divDnPn2}), it suffices to show that $ \s (H_i) = H_i$ for $i=1, \ldots , n$. For $i\neq j$, $\s (H_i-H_j) = H_i-H_j$, and so $d:= H_i'-H_i= H_j'-H_j$. For all $i=1, \ldots , n$,
$$ [\der_i , d] = \s ([\der_i, H_i]) - [\der_i, H_i]= \s (\der_i ) - \der_i = \der_i - \der_i =0.$$ So, $d\in C_{\divnc}(\CD_n ) = \CD_n$ (Lemma \ref{db11Mar13}.(1)) and $d= \sum_{i=1}^n \l_i\der_i$ for some $\l_i\in K$. The elements $H_1'=H_1+d, \ldots , H_n'=H_n+d$ commute hence $d=0$. Therefore, $\s = e$.  $\Box $

$\noindent $

{\bf Proof of Theorem \ref{A16Mar13}}. If $n=1$ then $\mathbf{G}_1^c\simeq \mT^1\ltimes \Sh_1$.  So, let $n\geq 2$. By (\ref{G=GF}) and Lemma \ref{a18Mar13}, $\mfGnc = G_n$. $\Box$

\begin{theorem}\label{A17Mar13}
For all $\s \in G_n$,
$$ \divv (\s (H_i )-H_i)=0\;\; {\rm for}\;\; i=1, \ldots , n.$$
\end{theorem}

{\it Proof}. For each $i=1, \ldots , n$, by (\ref{divDnPn2}), $\s (H_i) = \l_i (\s ) H_i +d_i(\s )$ for some elements $\l_i(\s ) \in K$ and $d_i(\s )\in \divn0$. By Theorem \ref{B16Mar13},
$$ 1=\divv (H_i) = \s (\divv (H_i)) = \divv (\s (H_i) ) = \divv (\l_i (\s ) H_i +d_i(\s )) = \l_i(\s ).\;\; \Box$$


{\bf The automorphisms of the Lie algebra $\divnc$ preserve divergence.}

\begin{corollary}\label{b17Mar13}
\begin{enumerate}
\item For all $\s \in \mfGnc$ and $\der \in \divnc$,
$ \divv (\s (\der ))= \divv (\der )$.
\item Every automorphism of the Lie algebra $\divn0$ is extendable to an automorphism of the Lie  algebra $\divnc$, and the extension is unique if $n\geq 2$.
\end{enumerate}
\end{corollary}

{\it Proof}. 1. The result follows from Theorem \ref{A16Mar13} and Theorem \ref{B16Mar13}.

2. Statement 2 follows from Theorem \ref{16Mar13} and  Theorem \ref{A16Mar13}. $\Box $


\begin{lemma}\label{b2Apr13}
For all $n\geq 2$, $P_n/K$ is a simple $\divn0$-module/$\divnc$-module; $\End_{\divn0}(P_n/K)=\End_{\divnc}(P_n/K)=K\id$ where $\id$ is the identity map. For $n=1$, $P_1/K$ is neither a  simple $\mathfrak{div}_1^0$-module nor a simple $\mathfrak{div}_1^c$-module.
\end{lemma}

{\it Proof}. If $n=1$ then $\mathfrak{div}_1^0=K\der_1$, $\mathfrak{div}_1^c=K\der_1+KH_1$ and $\sum_{i\leq m} Kx_1^i\der_1$ where $m\in \N$ are distinct $\mathfrak{div}_1^0/\mathfrak{div}_1^c$-submodules of $P_1$. Suppose that $n\geq 2$. It suffices to prove the statement for $\divn0$. Let $M$ be a nonzero $\divn0$-submodule of $P_n/K$ and $0\neq m\in M$. Using the actions of $\der_1, \ldots , \der_n\in \divn0$ on $m$ we obtain an element of $M$ of the form $\l x_i+K$ for some $\l\in K^*$. Hence, $x_i+K\in M$. Then $x_j+K= x_j\der_i *(x_i+K)\in M$ for all $j\neq i$. So, $x^\alpha +K\in M$ for all $\alpha \in \N^n$ with $|\alpha | =1$. We use induction on $|\alpha |$ to show that all $x^\alpha +K\in M$. Suppose that $m:=|\alpha | >1$.  If $x^\alpha +k = x_i^m +K$ for some $m\geq 2$ and $i$ then $x_i^m +K= (m-1)^{-1}  x_i(H_i-2H_j)*( x_i^{m-1}+K)$. Then, by applying elements of the type $x_j\der_i$ where $j\neq i$ to the element $x_i^m+K$ we obtain all the elements $x^\alpha +K$ with $|\alpha |=m$. Therefore, $P_n/K$ is a simple $\divn0$-module/$\divnc$-module.

Let $f\in \End_{\divn0}(P_n/K)$. Then applying $f$ to the equalities $\der_i*(x_1+K)=\d_{i1}$ for $i=1, \ldots , n$, we obtain the equalities
$$ \der_i*f(x_1+K)=\d_{i1} \;\; {\rm for }\;\; i=1, \ldots , n.$$ Hence, $f(x_1+K)\in \bigcap_{i=2}^n \ker_{P_n/K}(\der_i) \cap \ker_{P_n/K}(\der_1^2) = (K[x_1]/K)\cap \ker_{P_n/K}(\der_1^2) =K(x_1+K)$. So, $f(x_1+K) = \l ( x_1+K)$ and so $f=\l\,  \id$, by the simplicity of the $\divn0$-module  $P_n/K$. Therefore, $\End_{\divn0}(P_n/K)=\End_{\divnc}(P_n/K)=K\id$.
  $\Box $


\begin{proposition}\label{a2Apr13}
For $n\geq 2$, $\divnc$ is a maximal Lie subalgebra of $D_n$. For $n=1$, $\mathfrak{div}_1^c$ is not a maximal Lie subalgebra of $D_1$. For each $n\geq 1$, $\divnc$ is a $\mG_n$-invariant/$G_n$-invariant Lie subalgebra of $D_n$.
\end{proposition}

{\it Proof}. For $n=1$, $\mathfrak{div}_1^c=K\der_1+KH_1$ is contained in the Lie subalgebra $K\der_1+KH_1+Kx_1H_1$ of $D_1$. Suppose that $n\geq 2$. By (\ref{divDnPn}) and (\ref{divDnPn1}),
\begin{equation}\label{divDnPn3}
0\ra \divnc\ra D_n\stackrel{\divv}{\ra} P_n/ K\ra 0
\end{equation}
is the short exact sequence of $\divnc$-module. By Lemma \ref{b2Apr13}, the $\divnc$-module $P_n/K$ is simple. Then, $\divnc$ is a maximal Lie subalgebra of $D_n$.

By Theorem \ref{A16Mar13}, Theorem \ref{11Mar13} and Theorem \ref{B16Mar13}, $\divnc$ is a $\mG_n$-invariant/$G_n$-invariant Lie subalgebra of $D_n$.  $\Box $


\begin{lemma}\label{c2Apr13}
Let $n\geq 2$. Then
\begin{enumerate}
\item $P_n/K$ is a simple $G_n$-module with $\End_{G_n} (P_n/K)\simeq K$.
\item $D_n/ \divnc \stackrel{\divv}{\simeq}
P_n/ K$, an isomorphism of  $G_n$-modules.
\item $D_n/ \divnc$  is a simple $G_n$-module with $\End_{G_n}(P_n/K)\simeq K$.
\end{enumerate}
\end{lemma}

{\it Proof}. 1. Let $M$ be a $G_n$-submodule of $P_n$ properly containing $K$. We have to show that $M= P_n$, i.e. $x^\alpha \in M$ for all $\alpha \in \N^n$. The polynomial algebra $P_n=\oplus_{\alpha \in \N^n} Kx^\alpha$ is the direct sum of 1-dimensional non-isomorphic $\mT^n$-modules.  Hence $M$ is a homogeneous submodule of $P_n$. Hence $x^\beta \in M$ for some $\beta \neq 0$. If $\beta_i\neq 0$ then using the automorphism $s_i: x_i\mapsto x_i+1$, $x_j\mapsto x_j$, $j\neq i$, we see that $M\ni s_i(x^\beta ) - x^\beta$ and $\deg_{x_i} ( s_i(x^\beta ) - x^\beta)= \beta_i-1$. The module $M$ is closed under the maps $s_i-1$ for $i=1, \ldots , n$.  Hence, $x^\g\in M$ for all $\g \leq \beta$ where $\g \leq \beta$ iff $\g_i\leq \beta_i$ for all $i$. In particular, all $x_1, \ldots , x_n \in M$. Then applying the automorphism $\s_m :x_1\mapsto x_i+x_2^m$, $x_i\mapsto x_i$ for $i\neq 1$,  to the element $x_1$, we see that $M\ni (\s_m -1) (x_1) = x_2^m$ for all $m\geq 1$. Then applying the automorphism $ x_2\mapsto \sum_{i=1}^n x_i$, $x_i\mapsto x_i$ for $i\neq 2$, we have $(x_1+\cdots + x_n)^m\in M$. This implies that all $x^\alpha \in M$, by the homogeneity of $M$.

Let $f\in \End_{G_n}(P_n/K)$. Since $f$ commutes with the action of the subgroup $\mT^n$ of $G_n$, we must have $ f(x^\alpha +K) = \l_\alpha (x^\alpha +K)$ for all $\alpha \in \N^n$ and some $\l_\alpha \in K$. In particular, $f(x_1+K) = \l (x_1+K)$ for some $\l \in K$. Since $f$ commutes with the action of the symmetric group $S_n$ (which is obviously a subgroup of $G_n$), $f(x_i+K)= \l (x_i+K)$ for all $i=1, \ldots , n$. Now, we use induction on $|\alpha |$ show that $f(x^\alpha +K) = \l ( x^\alpha +K)$. The initial case when $|\alpha |=1$ has just been established. So, let $|\alpha |>1$. Then $\alpha_i>0$ for some $i$, and $\deg ((s_i-1)x^\alpha +K) <|\alpha |$. By induction, $f((s_i-1) x^\alpha +K) = \l ((s_i-1) x^\alpha +K)$. Now, it follows from the equality
$$ f(s_i(x^\alpha ) +K) -f( x^\alpha +K) = f((s_i-1)x^\alpha +K) = \l (s_i-1)x^\alpha +K$$ that $\l_\alpha = \l$, and so $f= \l\,  \id$. Therefore, $\End_{G_n}(P_n/K)= K\id$.

2. Statement 2 follows from (\ref{divDnPn3}) and Theorem \ref{B16Mar13}.

3. Statement 3 follows from statements 1 and 2. $\Box $

\begin{lemma}\label{a14Apr13}
The Lie algebra $\divn0$ is a $\mG_n$-invariant/$G_n$-invariant Lie subalgebra of $D_n$.
\end{lemma}

{\it Proof}. The statement follows from Theorem \ref{16Mar13}, Theorem \ref{11Mar13} and Theorem \ref{B16Mar13}. $\Box $

$\noindent $

{\bf Conjecture}: {\it Every nonzero homomorphism of the Lie algebra $\divn0$ is an automorphism.}



\section{Minimal set of generators for the Lie algebras $\divn0$ and $\divnc$}\label{TMMSG}

 In this section, the proofs of Theorem \ref{6Oct13} and Theorem \ref{A6Oct13} are given.

\begin{equation}\label{ditai2}
[x_i^2\der_{i+1}, x_{i+1}\der_i]=\th_i^{e_i}.
\end{equation}
In more detail, LHS$= x_i^2\der_i-2x_ix_{i+1}\der_{i+1}=x_i(H_i-2H_{i+1}) = \th_i^{e_i}$.
\begin{equation}\label{ditai3}
[x_{i+1}^2\der_i, x_i\der_{i+1}]=-\th_i^{e_{i+1}}.
\end{equation}
Similarly, LHS$= x_{i+1}^2\der_{i+1}-2x_ix_{i+1}\der_i=-x_{i+1}(2H_i-H_{i+1}) = -\th_i^{e_{i+1}}$.

$\noindent $

{\bf Proof of Theorem \ref{6Oct13}}. The elements in Theorem \ref{6Oct13} belong to $\divn0$ and let $\CG$ be the Lie subalgebra of $\divn0$ they generate.

(i) $\CG = \divn0$: To prove that the equality holds we use induction on $n$. Let $n=2$. So, $ \CG = \langle \der_1, x_2^2\der_1, x_1^2\der_2\rangle$. Then $\der_2, x_1\der_2, x_2\der_1\in \CG$:
$$ x_1\der_2= \frac{1}{2}[\der_1, x_1^2\der_2], \;\; \der_2= [ \der_1, x_1\der_2] , \;\; x_2\der_1= \frac{1}{2}[\der_2, x_2^2\der_1].$$
By Lemma \ref{a17Mar13}.(3), we have to show that the elements $x_2^i\der_1$, $x_1^i \der_2$ and $\th_1^\alpha = x^\alpha ((\alpha_2+1)H_2-(\alpha_1+1)H_2)$ belong to $\CG$ where $i\in \N$ and $\alpha \in \N^2$. By (\ref{ditai2}) and (\ref{ditai3}),
$$\th_1^{e_1} = [ x_1^2\der_2, x_2\der_1] \in \CG , \;\; \th_1^{e_2} = -[x_2^2\der_1, x_1\der_2] \in \CG.$$Then using (\ref{ditai}) and (\ref{ditai1}), we see that $\th_1^\alpha \in \CG$ for all $0\neq \alpha \in \N^2$. For $\alpha =0$, $\th_1^0 = H_1-H_2= [ x_1\der_2, x_2\der_1] \in \CG$. By (\ref{xadxbd2}), for $j\geq 1$,
$$ \CG \ni [x_2\der_1, \th_1^{je_2}]= -((j+1)H_1-H_2, -e_1+e_2)x_2^{j+1}\der_1= (j+2) x_2^{j+1}\der_1.$$ By symmetry, $\CG \ni (j+2) x_1^{j+1}\der_2$, i.e. $\CG = \divn0$.

Suppose that $n\geq 3$, and that the equality $\CG =\mathfrak{div}_{n'}^0 $ holds for all $n'$ such that $2\leq n' <n$.

{\em Step 1}. $\{ \der_i, x_j^\nu \der_k\, | \, i,j,k=1, \ldots , n; j\neq k; \nu =1,2\} \subseteq \CG$: For $i=2, \ldots , n$, $\der_i = \frac{1}{2}[\der_1, [\der_1, x_1^2\der_i]]\in \CG$ and $x_1\der_i = \frac{1}{2} [\der_1, x_1^2\der_i]\in \CG$. Then $x_i\der_1= \frac{1}{2} [\der_i,x_i^2\der_1]\in \CG$ for $i=2, \ldots , n$. For all $i\neq j$ such that $i\neq 1$ and $j\neq 1$, $x_i\der_j =[x_i\der_1, x_1\der_j]\in \CG$ and $ x_i^2\der_j = [x_i^2\der_1, x_1\der_j]\in \CG$. For $i=2, \ldots , n$, fix an index $j$ such that $j\neq 1, i$. Then $x_i\der_1= [ x_i\der_j , x_j \der_1]\in \CG$ and $ x_i^2\der_1= [ x_i^2\der_j , x_j \der_1] \in \CG$. The proof of Step 1 is complete.

For $i=2,\ldots , n$, let $\mathfrak{div}^0_{n,i}$ be the Lie algebra $\mathfrak{div}^0_{n-1}$ for the polynomial algebra $K[x_1, \ldots , \widehat{x_i}, \ldots , x_n]$ ($x_i$ is missed).

{\em Step 2}. {\em   For} $i=1, \ldots , n$, $\mathfrak{div}^0_{n,i}\subseteq \CG$: This follows from Step 1 and induction.

{\em Step 3}. $\din0\subseteq\CG$: This follows from (\ref{xadxbd1}) and Step 2.

{\em Step 4}. $\ivn0\subseteq\CG$: This follows from (\ref{xadxbd2}) (where $\alpha_i=0$) and Step 3.

(ii) {\em Minimality}:

{\em (a) The element $\der_1$ cannot be dropped}: By Lemma \ref{a17Mar13}.(3), $\divn0$ is a $\Z$-graded Lie algebra which is determined by the degree $\deg$ in the following way,  for $x^\alpha\der_i\in P_n^{\der_i} \der_i$, $\deg (x^\alpha\der_i)= |\alpha |-1$ and $\deg (\th_i^\alpha ) = |\alpha |$ for all $\alpha \in \N^n$ and $i=1, \ldots , n-1$. Clearly, $\deg (\der_1) = -1$ (negative) and the degrees the rest of the generators are equal to 1 (positive). Therefore, the $\der_1$ cannot be dropped.

{\em (b) The element $x_i^2\der_1$ ($i=2, \ldots , n$) cannot be dropped}: Since otherwise the Lie algebra generated by the remaining elements would belong to $\oplus_{j=1}^n K [x_1, \ldots , \widehat{x_i}, \ldots , x_n]\der_j$ ($x_i$ is missed), a contradiction (see (i)).

{\em (c) The element $x_1^2\der_i$ ($i=2, \ldots , n$) cannot be dropped}: Since otherwise the Lie algebra generated by the remaining elements would belong to $\oplus_{j\neq i}^n P_n\der_j$, a contradiction (see (i)). $\Box $

$\noindent $

{\bf Proof of Theorem \ref{A6Oct13}}. By Theorem \ref{6Oct13}, the elements in Theorem \ref{A6Oct13} generate the Lie algebra $\divnc = \divn0 \oplus KH_1$ and the element $H_1$ cannot be dropped (by Theorem \ref{6Oct13}).

{\em (a) The element $\der_1$ cannot be dropped} by the same reason as in the proof of Theorem \ref{6Oct13} as $\deg (H_1)=0$.

{\em (b) The element $x_i^2\der_1$ ($i=2, \ldots , n$) cannot be dropped}  by the same reason as in the proof of Theorem \ref{6Oct13}.

{\em (c) The element $x_1^2\der_i$ ($i=2, \ldots , n$) cannot be dropped}  by the same reason as in the proof of Theorem \ref{6Oct13}. $\Box $


$${\bf Acknowledgements}$$

 The work is partly supported by  the Royal Society  and EPSRC.

\small{

Department of Pure Mathematics

University of Sheffield

Hicks Building

Sheffield S3 7RH

UK

email: v.bavula@sheffield.ac.uk}

\end{document}